\documentclass{amsart}%
\usepackage{amsmath}
\usepackage{graphicx}
\usepackage{amsfonts}
\usepackage{amssymb}%
\newtheorem{theorem}{Theorem}[section]
\theoremstyle{plain}

\newtheorem{corollary}[theorem]{Corollary}

\newtheorem{lemma}[theorem]{Lemma}

\newtheorem{proposition}[theorem]{Proposition}
\newtheorem{remark}[theorem]{Remark}
\numberwithin{equation}{section}

\begin{document}
\title[The Miura Map on the Line]{The Miura Map on the Line}
\author{Thomas Kappeler}
\address[Kappeler]{ Department of Mathematics, University of Z\"{u}rich, Z\"{u}rich, Switzerland}
\email{tk@math.unizh.ch}
\author{Peter Perry}
\address[Perry]{ Department of Mathematics, University of Kentucky, Lexington, Kentucky
40506-0027, U. S. A.}
\email{perry@ms.uky.edu}
\author{Mikhail Shubin}
\address[Shubin]{ Department of Mathematics, Northeastern University, Boston,
Massachusetts 02115, U. S. A.}
\email{shubin@neu.edu}
\author{Peter Topalov}
\address[Topalov]{ Institute of Mathematics, Bulgarian Academy of Sciences, Sofia, Bulgaria.}
\curraddr{Department of Mathematics, University of Z\"{u}rich, Z\"{u}rich, Switzerland}
\email{topalov@math.bas.bg, topalov@math.unizh.ch}
\thanks{Thomas Kappeler supported in part by the Swiss National Science Foundation. }
\thanks{Peter Perry supported in part by NSF grants DMS-0100829 and DMS-0408419.}
\thanks{Mikhail Shubin supported in part by NSF\ grant DMS-0107796.}
\thanks{Peter Topalov supported in part by the European Training Network HPRN-CT-00118.}

\begin{abstract}
We study relations between properties of the Miura map $r\mapsto
q=B(r)=r^{\prime}+r^{2}$ and Schr\"{o}dinger operators $L_{q}=-d^{2}/dx^{2}+q$
where $r$ and $q$ are real-valued functions or distributions (possibly not
decaying at infinity) from various classes. In particular, we study $B$ as a
map from $L_{\mathrm{loc}}^{2}(\mathbb{R})$ to the local Sobolev space
$H_{\mathrm{loc}}^{-1}(\mathbb{R})$ and the restriction of $B$ to the Sobolev
spaces $H^{\beta}(\mathbb{R})$ with $\beta\geq0$. For example, we prove that
the image of $B$ on $L_{\mathrm{loc}}^{2}(\mathbb{R})$ consists exactly of
those $q\in H_{\mathrm{loc}}^{-1}(\mathbb{R})$ such that the operator $L_{q}$
is positive. We also investigate mapping properties of the Miura map in these
spaces. As an application we prove an existence result for solutions of the
Korteweg-de Vries equation in $H^{-1}(\mathbb{R})$ for initial data in the
range $B(L^{2}(\mathbb{R}))$ of the Miura map.

\end{abstract}
\maketitle
\tableofcontents

\section{Introduction and Main Results}

\label{sec.intro}

The Miura map is the nonlinear mapping
\begin{equation}
B(r)=r^{\prime}+r^{2} \label{eq.B}%
\end{equation}
which takes classical solutions of the modified Korteweg-de Vries (mKdV)
equation to classical solutions of the Korteweg-de Vries (KdV)\ equation. More
precisely, let%
\[
\operatorname{mKdV}(v):=v_{t}-6v^{2}v_{x}+v_{xxx}%
\]
and%
\[
\operatorname*{KdV}(u):=u_{t}-6uu_{x}+u_{xxx}%
\]
for functions $v,u\in\mathcal{C}^{\infty}(\mathbb{R}\times\mathbb{R})$. Then
\begin{equation}
\operatorname*{KdV}(B(v))=\left(  \operatorname{mKdV}(v)\right)  _{x}%
+2v\cdot\operatorname{mKdV}(v) \label{eq.MiuraId}%
\end{equation}
so that $\operatorname*{KdV}(B(v))=0$ whenever $mKdV(v)=0$. The Miura map has
been used extensively to relate existence and uniqueness results for solutions
for the mKdV and KdV equations. More fundamentally, the global geometry of
nonlinear differentiable mappings such as the Miura map has been studied to
solve various nonlinear differential equations.

It is well-known that the KdV equation can be successfully studied with the
help of the spectral theory of the Schr\"{o}dinger operators
\begin{equation}
L_{q}:=-d^{2}/dx^{2}+q. \label{E:Schroedinger}%
\end{equation}
In particular, in appropriate classes of potentials $q$, the spectrum of
$L_{q}$ is preserved by the KdV flow applied to $q$. (See e.g. \cite{Lax:1968}%
.) The Miura map is formally related with the Schr\"{o}dinger operator as
follows: the relation $q=B(r)=r^{\prime}+r^{2}$ is equivalent to the following
factorization of $L_{q}$:
\begin{equation}
L_{q}=L_{B(r)}=(\partial_{x}-r)^{+}(\partial_{x}-r)=(-\partial_{x}%
-r)(\partial_{x}-r), \label{E:factor}%
\end{equation}
where $\partial_{x}=d/dx$, and $A^{+}$ means the operator formally adjoint to
an operator $A$ in functions on $\mathbb{R}$ (with respect to the scalar
product in $L^{2}({\mathbb{R}})$).

The aim of this paper is to study the Miura map and its geometry on the real
line with an emphasis on function spaces of low regularity. Our first result
characterizes the range of the Miura map. We denote by $B$ the map
(\ref{eq.B}) from real-valued functions in $L_{\mathrm{loc}}^{2}(\mathbb{R})$
into the local Sobolev space $H_{\mathrm{loc}}^{-1}(\mathbb{R})$. If $q$ is a
real-valued distribution in $H_{\mathrm{loc}}^{-1}(\mathbb{R})$, the operator
$L_{q}$ maps $\mathcal{C}_{0}^{\infty}(\mathbb{R})$ into the Sobolev space
$H^{-1}(\mathbb{R})$, so that $\left(  L_{q}\varphi,\varphi\right)  $ is
well-defined for any $\varphi\in\mathcal{C}_{0}^{\infty}(\mathbb{R})$. We
write $L_{q}\geq0$ if $\left(  L_{q}\varphi,\varphi\right)  \geq0$ for all
$\varphi\in\mathcal{C}_{0}^{\infty}(\mathbb{R})$.

\begin{theorem}
\label{thm.im1}Let $q$ be a real-valued distribution belonging to
$H_{\mathrm{loc}}^{-1}(\mathbb{R})$. Then the following statements are
equivalent. \newline(i) $q\in\operatorname{Im}(B)$, i.e., $q=r^{\prime}+r^{2}$
for some real-valued function $r\in L_{\mathrm{loc}}^{2}(\mathbb{R})$.
\newline(ii) The equation $L_{q}y=0$ has a strictly positive solution $y\in
H_{\mathrm{loc}}^{1}(\mathbb{R})$. \newline(iii) $L_{q}\geq0$.
\end{theorem}

Next, we consider the restriction of $B$ to the Sobolev space $H^{\beta
}(\mathbb{R})$ for $\beta\geq0$; we denote this restriction by $B_{\beta}$.
Although $B_{\beta}$ has range contained in $H^{\beta-1}(\mathbb{R})$, it is
not true that $\operatorname{Im}(B_{\beta})=\operatorname{Im}(B)\cap
H^{\beta-1}(\mathbb{R})$; rather, an additional condition is needed to
characterize its range.

\begin{theorem}
\label{thm.im2}Let $\beta\ge0$ be arbitrary. A real-valued distribution $q\in
H^{\beta-1}(\mathbb{R})$ belongs to $\operatorname{Im}(B_{\beta})$ if and only
if \newline(i) $L_{q}\geq0$, and\newline(ii) $q$ can be presented in the form
$q=f^{\prime}+g$ for $f\in L^{2}(\mathbb{R})$ and $g\in L^{1}\left(
\mathbb{R}\right)  $.
\end{theorem}

In addition, we will give an alternative characterization of
$\operatorname{Im}(B_{\beta})$ in terms of a \textquotedblleft special
integral\textquotedblright\ of $q$ on $\mathbb{R}$ which coincides with the
ordinary integral of $q$ if $q\in L^{1}\left(  \mathbb{R}\right)  $ (see
Theorem \ref{thm.special}). 

We also study the geometry of the Miura map. Whereas the Miura map
on spaces of periodic functions is known to be a global fold, the situation
for non-periodic functions is completely different. Let%
\begin{equation}
E_{1}=\left\{  q\in\operatorname{Im}(B):B^{-1}(q)\text{ consists of a single
point}\right\}  \label{eq.e1}%
\end{equation}
and%
\begin{equation}
E_{2}=\left\{  q\in\operatorname{Im}(B):B^{-1}(q)\text{ is homeomorphic to an
interval}\right\}  . \label{eq.e2}%
\end{equation}
Here we consider $B^{-1}(q)$ with the natural Fr\'echet topology of
$L_{\mathrm{loc}}^{2}(\mathbb{R})$.

\begin{theorem}
\label{thm.im3}$\operatorname{Im}(B)=E_{1}\cup E_{2}$ and both $E_{1}$ and
$E_{2}$ are dense in $\operatorname{Im}(B)$ in the natural Fr\'echet topology
of $H_{\mathrm{loc}}^{-1}(\mathbb{R})$.
\end{theorem}

As an application of our results on the Miura map, we prove an existence
result for solutions of the Korteweg-de Vries equation in $H^{-1}(\mathbb{R})$
for initial data in the range $\operatorname{Im}(B_{0})$ of the Miura map
$B_{0}:L^{2}(\mathbb{R})\rightarrow H^{-1}(\mathbb{R})$. We follow the
approach of Tsutsumi \cite{Tsu:1989}, \ who proved such an existence result
when the initial data is a positive, finite Radon measure on $\mathbb{R}$. His
arguments, combined with our results on the Miura map $B_{0},$ lead to the
following theorem.

\begin{theorem}
\label{thm.KdV4}Assume that $u_{0}\in\operatorname{Im}(B_{0})$. Then there
exists a global weak solution of KdV with $u(t)\in\operatorname{Im}\left(
B_{0}\right)  $ for all $t\in\mathbb{R}$. More precisely: \newline(i) $u\in
L^{\infty}(\mathbb{R},H^{-1}(\mathbb{R}))\cap L_{\mathrm{loc}}^{2}%
(\mathbb{R}^{2})$; \newline(ii) for all functions $\varphi\in\mathcal{C}%
_{0}^{\infty}(\mathbb{R}^{2})$, the identity%
\[
\int_{\mathbb{R}}\int_{\mathbb{R}}\left(  -u\varphi_{t} - u\varphi_{xxx}%
+3u^{2}\varphi_{x}\right)  ~dx~dt=0
\]
holds, and \newline(iii) $\lim_{t\rightarrow0}u(t)=u_{0}$ in $H^{-1}%
(\mathbb{R})$.
\end{theorem}

We state and prove a slightly stronger version of the above theorem in Section
\ref{sec.KdV}.

In Appendix \ref{app.work}, we provide a few comments on the Miura map as well
as on other work related to the results presented in this paper.

\medskip

\noindent\textbf{Acknowledgements}. Peter Perry and Mikhail Shubin thank the
Mathematical Institute at the University of Z\"{u}rich for its hospitality
during part of the time that this work was done. The authors thank Carlos
Tomei for bringing the paper of McKean and Scovel \cite{MS:1986} to their
attention. They also thank Vladimir Maz'ya and Tatyana Shaposhnikova for
pointing out the reference \cite{RS:1996} (see the proof of Lemma
\ref{lemma.pr}).

\section{Preliminaries}

\subsection{Spaces of Distributions}

For $\alpha\in\mathbb{R}$, we denote by $H^{\alpha}(\mathbb{R})$ the
completion of $\mathcal{C}_{0}^{\infty}(\mathbb{R})$ in the norm%
\[
\left\Vert \varphi\right\Vert _{\alpha}=\left(  \int\left(  1+\left\vert
\xi\right\vert ^{2}\right)  ^{\alpha}\left\vert \widehat{u}(\xi)\right\vert
^{2}~\frac{d\xi}{2\pi}\right)  ^{1/2}%
\]
where%
\[
\widehat{u}(\xi)=\int\exp\left(  -i\xi x\right)  ~u(x)~dx.
\]
Clearly, $H^{0}(\mathbb{R})=L^{2}(\mathbb{R})$ and $H^{\alpha}(\mathbb{R}%
)\subset H^{\beta}(\mathbb{R})$ if $\alpha\geq\beta$. Therefore $H^{\alpha
}(\mathbb{R})\subset L^{2}(\mathbb{R})$ if $\alpha\geq0$. For any $\alpha
\in\mathbb{R}$, $H^{\alpha}(\mathbb{R})$ is a space of tempered distributions.
A distribution $u$ belongs to $H_{\mathrm{loc}}^{\alpha}(\mathbb{R})$ if $\chi
u\in H^{\alpha}(\mathbb{R})$ for any function $\chi\in\mathcal{C}_{0}^{\infty
}(\mathbb{R})$. We consider $H_{\mathrm{loc}}^{\alpha}(\mathbb{R})$ and
$L_{\mathrm{loc}}^{p}(\mathbb{R})$ with $p\geq1$ with their natural
Fr\'{e}chet topology - see \cite{Treves:1967}, chapter 31-12.


A classical result of distribution theory (see, for example, chapter 1 of
\cite{GS:1964}) asserts that if $q\in H_{\mathrm{loc}}^{-1}(\mathbb{R})$, then
$q=Q^{\prime}$ for a function $Q\in L_{\mathrm{loc}}^{2}(\mathbb{R})$. If
$q\in H^{\beta-1}(\mathbb{R})$ for $\beta\geq0$ we have a sharper result. Let%
\[
H^{\infty}(\mathbb{R})=%
{\displaystyle\bigcap\limits_{\beta\geq0}}
H^{\beta}(\mathbb{R})\subset\mathcal{C}^{\infty}(\mathbb{R})
\]

\begin{lemma}
\label{lemma.rep}Let $\beta\geq0$ and let $q\in H^{\beta-1}(\mathbb{R})$.
There exist functions $f\in H^{\beta}(\mathbb{R})$ and $g\in H^{\infty
}(\mathbb{R})$ so that $q=f^{\prime}+g$ as elements of $H^{\beta-1}%
(\mathbb{R})$.
\end{lemma}

\begin{proof}
Let $\psi\in\mathcal{C}_{0}^{\infty}(\mathbb{R})$ with $\psi(\xi)=1$ near
$\xi=0$, and choose%
\[
\widehat{f}(\xi)=\left(  i\xi\right)  ^{-1}\left(  1-\psi(\xi)\right)
~\widehat{q}(\xi)
\]
and%
\[
\widehat{g}(\xi)=\psi(\xi)~\widehat{q}(\xi).
\]

\end{proof}

Finally, the following result will be useful in studying the continuity of the
Miura map and the regularity of solutions to the Riccati equation
$q=r^{\prime}+r^{2}$.

\begin{lemma}
\label{lemma.pr}The multiplication $\left\{  u,v\right\}  \mapsto uv$ can be
extended from the bilinear map $\mathcal{C}_{0}^{\infty}(\mathbb{R}%
)\times\mathcal{C}_{0}^{\infty}(\mathbb{R})\rightarrow\mathcal{C}_{0}^{\infty
}(\mathbb{R})$ to continuous bilinear maps%
\begin{equation}
H^{\beta}(\mathbb{R})\times H^{\beta}(\mathbb{R})\rightarrow H^{\beta
}(\mathbb{R}),~~\beta>1/2, \label{eq.H1}%
\end{equation}%
\begin{equation}
H^{1/2}(\mathbb{R})\times H^{1/2}(\mathbb{R})\rightarrow H^{1/2-\delta
}(\mathbb{R})\text{ for any }\delta>0, \label{eq.H2}%
\end{equation}%
\begin{equation}
H^{\beta}(\mathbb{R})\times H^{\beta}(\mathbb{R})\rightarrow H^{2\beta
-1/2}(\mathbb{R}),~~0<\beta<1/2, \label{eq.H3}%
\end{equation}%
\begin{equation}
L^{2}(\mathbb{R})\times L^{2}(\mathbb{R})\rightarrow L^{1}(\mathbb{R})\subset
H^{-1/2-\delta}(\mathbb{R)}\text{ for any }\delta>0\text{.} \label{eq.H4}%
\end{equation}

\end{lemma}

\begin{proof}
All of the statements, except the last one, are particular cases of more
general multidimensional results formulated, for example, in Theorem 1 of
section 4.6.1 of \cite{RS:1996}. The last statement is obvious except the last
inclusion; this follows by duality from the imbedding $H^{1/2+\delta
}(\mathbb{R})\rightarrow L^{\infty}(\mathbb{R})$, which is a particular case
of a well-known Sobolev imbedding theorem.
\end{proof}

Let us also introduce a notation $(\cdot,\cdot)$ for miscellaneous
sesquilinear pairings extending the pairing
\[
(u,v)=\int_{{\mathbb{R}}}u(x)\overline{v(x)} dx,\quad u,v\in\mathcal{C}%
_{0}^{\infty}({\mathbb{R}}),
\]
by continuity. In particular, we will use the extended pairings in the
following cases:

(i) $u$ is a distribution on ${\mathbb{R}}$, $v\in\mathcal{C}_{0}^{\infty
}({\mathbb{R}})$;

(ii) $u\in H^{s}({\mathbb{R}})$, $v\in H^{-s}({\mathbb{R}})$, where
$s\in{\mathbb{R}}$;

(iii) $u\in H^{s}_{\mathrm{l}oc}({\mathbb{R}})$, $v\in H^{-s}_{\mathrm{c}%
omp}({\mathbb{R}})$, where $s\in{\mathbb{R}}$ and $H^{-s}_{comp}({\mathbb{R}%
})$ is the space of compactly supported distributions from $H^{-s}%
({\mathbb{R}})$.

The integration by parts formula
\begin{equation}
\label{E:parts}(U^{\prime},v)=-(U,v^{\prime})
\end{equation}
holds in all these cases, e.g. for $U\in H^{s+1}_{\mathrm{l}oc}({\mathbb{R}}%
)$, $v\in H^{-s}_{\mathrm{c}omp}({\mathbb{R}})$, $s\in{\mathbb{R}}$, by
continuity of the pairing.

\subsection{Continuity of the Miura Map}

It is easy to see that the Miura map defines a bounded continuous map from
$L^{2}(\mathbb{R})$ into $H^{-1}(\mathbb{R})$ (we use the standard embedding
$L^{1}(\mathbb{R})\subset H^{-1}(\mathbb{R})$ which follows by duality from
the Sobolev embedding $H^{1}(\mathbb{R})\subset L^{\infty}(\mathbb{R})$).
Localizing this, we see that the Miura map defines a bounded continuous
mapping from $L_{\mathrm{loc}}^{2}(\mathbb{R})$ into $H_{\mathrm{loc}}%
^{-1}(\mathbb{R})$. This is extended to more general Sobolev spaces as follows.

\begin{proposition}
\label{prop.Miura.continuous}Let $\beta\geq0$. The Miura map is a continuous
mapping from $H^{\beta}(\mathbb{R})$ into $H^{\beta-1}(\mathbb{R})$ and from
$H_{\mathrm{loc}}^{\beta}(\mathbb{R})$ into $H_{\mathrm{loc}}^{\beta
-1}(\mathbb{R})$ which is bounded, i.e. maps bounded subsets into bounded subsets.
\end{proposition}

\begin{remark}
For the definition of bounded sets in Fr\'{e}chet spaces see e.g.
\cite{Treves:1967}, Chapter 14, especially Proposition 14.5.
\end{remark}

\begin{proof}
It suffices to prove the first statement since the second follows by
localization. It is clear that the map $r\mapsto r^{\prime}$ is a bounded
continuous map from $H^{\beta}(\mathbb{R})$ into $H^{\beta-1}(\mathbb{R})$ so
it suffices to show that the map $r\mapsto r^{2}$ is a bounded continuous map
from $H^{\beta}(\mathbb{R})$ into $H^{\beta-1}(\mathbb{R})$. This follows from
Lemma \ref{lemma.pr} and the trivial embedding of $H^{\alpha}(\mathbb{R})$
into $H^{\gamma}(\mathbb{R})$ if $\alpha\geq\gamma$.
\end{proof}

\subsection{A First-Order System}

Let $y\in H_{\mathrm{loc}}^{1}(\mathbb{R})$ be a solution of the equation
\begin{equation}
-y^{\prime\prime}+qy=0 \label{eq.sec}%
\end{equation}
where $q\in H_{\mathrm{loc}}^{-1}(\mathbb{R})$. Then we can conveniently
rewrite the equation in the form of a first-order system%
\begin{equation}
\left\{
\begin{array}
[c]{ccl}%
y^{\prime} & = & Qy + u\\
u^{\prime} & = & -Q^{2}y-Qu
\end{array}
\right.  \label{eq.sys}%
\end{equation}
where $Q\in L_{\mathrm{loc}}^{2}(\mathbb{R})$ is a real-valued function with
$Q^{\prime}=q$ (this is a well-known procedure in the study of differential
operators with singular coefficients; see \S 1.1 of \cite{SS:2003} and
references therein). This system is equivalent to equation (\ref{eq.sec}) in
the following sense. If $y\in H_{\mathrm{loc}}^{1}(\mathbb{R})$ is a solution
of (\ref{eq.sec}), then taking $u=y^{\prime}-Qy$ we obtain by straightforward
substitution that the pair $\left\{  y,u\right\}  $ satisfies (\ref{eq.sys})
in the sense of distributions. It follows that $u\in W_{\mathrm{loc}}%
^{1,1}(\mathbb{R})$, the space of $L_{\mathrm{loc}}^{1}(\mathbb{R})$-functions
with distributional derivatives in $L_{\mathrm{loc}}^{1}(\mathbb{R})$. On the
other hand, if $y$ and $u$ belong to $W_{\mathrm{loc}}^{1,1}(\mathbb{R})$ and
the pair $\left\{  y,u\right\}  $ satisfies (\ref{eq.sys}), then in fact $y\in
H_{\mathrm{loc}}^{1}(\mathbb{R})$ and $y$ satisfies (\ref{eq.sec}).

Since the coefficients of the linear system (\ref{eq.sys}) are in
$L_{\mathrm{loc}}^{1}(\mathbb{R})$, the standard existence and uniqueness
result holds for the corresponding initial value problem on the whole real
line, and the solutions $\left\{  y,u\right\}  $ depend continuously on the
initial data.

The next lemma shows that nonnegative solutions of $L_{q}y=0$ are either
strictly positive or identically zero.

\begin{lemma}
\label{lemma.zero} Suppose that $y\in H_{\mathrm{loc}}^{1}(I)$ is a solution
of $L_{q}y=0$ with a real-valued $q\in H_{\mathrm{loc}}^{-1}(I)$, where $I$ is
an open interval in ${\mathbb{R}}$. Assume that $y(x_{0})=0$ for some
$x_{0}\in I$. Denote $u=y^{\prime}-Qy$ as in \eqref{eq.sys}. Then the
following statements hold true:

A. We have the following trichotomy for the behavior of $y$ near $x_{0}$:

(i) If $u(x_{0})=0$, then $y\equiv0$ on $I$.

(ii) If $u(x_{0})>0$ then in a neighborhood of $x_{0}$ we have $y(x)<0$ for
$x<x_{0}$, and $y(x_{0})>0$ for $x>x_{0}$.

(iii) If $u(x_{0})<0$ then in a neighborhood of $x_{0}$ we have $y(x)>0$ for
$x<x_{0}$, and $y(x_{0})<0$ for $x>x_{0}$.

B. If $y\not \equiv 0$, then all zeros of $y$ on $I$ are isolated.

C. If $y(x)\geq0$ near $x_{0}$ then $y(x)\equiv0$ on $I$.
\end{lemma}

\begin{proof}
Clearly, $(i)$ follows from the uniqueness of solution $\{y,u\}$ of
\eqref{eq.sys} with the given initial conditions $y(x_{0})=y_{0},
u(x_{0})=u_{0}$.

Define $z(x)= y(x)\exp\left(  -\int_{x_{0}}^{x}Q(s)~ds\right)  $ and note
that, as $u$ is absolutely continuous, so is $z^{\prime}(x)=u(x)\exp\left(
-\int_{x_{0}}^{x}Q(s)~ds\right)  $. Clearly, $z(x_{0})=0$, and $z(x)$ has the
same sign as $y(x)$ for all $x\in I$. Also, $u(x_{0})>0$ is equivalent to
$z^{\prime}(x_{0})>0$, so $(ii)$ and $(iii)$ immediately follow.

Clearly, $B$ and $C$ follow from $A$.
\end{proof}

\section{Positivity, Positive Solutions, and the Miura Map}

In this section we prove Theorem \ref{thm.im1}. First, we describe the
connection between the Miura map and positive solutions of $L_{q}y=0$. For a
real-valued distribution $q\in H_{\mathrm{loc}}^{-1}(\mathbb{R})$, let
$\operatorname*{Pos}(q)$ denote the (possibly empty) set of functions $y\in
H_{\mathrm{loc}}^{1}(\mathbb{R})$ with the properties that $L_{q}y=0$,
$y(x)>0$ for all $x\in\mathbb{R}$, $\ $and $y(0)=1$.

\begin{lemma}
\label{lemma.straight}Let $q$ be a real-valued distribution belonging to
$H_{\mathrm{loc}}^{-1}(\mathbb{R})$.\newline(a) If $y\in\operatorname*{Pos}%
(q)$ then $q=B(y^{\prime}/y)$.\newline(b) If $q=B(r)$ for some $r\in
L_{\mathrm{loc}}^{2}(\mathbb{R})$ then $y(x)=\exp\left(  \int_{0}%
^{x}r(s)~ds\right)  $ belongs to $\operatorname*{Pos}(q)$.
\end{lemma}

The proof is easily obtained by straightforward calculations.

The maps%
\[
\operatorname*{Pos}(q)\ni y\mapsto\frac{d}{dx}\log(y(x))\in B^{-1}(q)
\]
and
\[
B^{-1}(q)\ni r\mapsto\exp\left(  \int_{0}^{x}r(s)~ds\right)  \in
\operatorname*{Pos}(q)
\]
are continuous if we topologize $B^{-1}(q)$ with the topology induced from
$L_{\mathrm{loc}}^{2}(\mathbb{R})$ and $\operatorname*{Pos}(q)$ with the
topology induced from $H_{\mathrm{loc}}^{1}(\mathbb{R})$. These maps are
mutual inverses. Hence, we have shown:

\begin{proposition}
\label{prop.homeo}The set $B^{-1}(q)$ is nonempty if and only if
$\operatorname*{Pos}(q)$ is nonempty. For any $r\in B^{-1}(q)$, $B^{-1}(B(r))$
is homeomorphic to $\operatorname*{Pos}(B(r))$.
\end{proposition}

Next, we show that if $L_{q}\geq0$, then $\operatorname*{Pos}(q)$ is nonempty.
To this end, we introduce the sesquilinear forms%
\begin{equation}
\label{E:tq}\mathfrak{t}_{q}(\varphi,\psi)=\int_{\mathbb{R}}\varphi^{\prime
}(x)\overline{\psi^{\prime}}(x)dx+\left(  q,\overline{\varphi}\psi\right)
\end{equation}
and%
\begin{equation}
\label{E:tqI}\mathfrak{t}_{q,I}(\varphi,\psi)=\int_{I}\varphi^{\prime
}(x)\overline{\psi^{\prime}}(x)~dx+\left(  q,\overline{\varphi}\psi\right)
\end{equation}
defined respectively on $\mathcal{C}_{0}^{\infty}(\mathbb{R})$ and
$\mathcal{C}_{0}^{\infty}(I)$, where $I=(a,b)$ is a bounded, open interval of
$\mathbb{R}$.

The form $\mathfrak{t}_{q}$ is also well defined by \eqref{E:tq} for
$\varphi,\psi\in H^{1}_{comp}({\mathbb{R}})$, where $H^{1}_{comp}({\mathbb{R}%
})$ is the space of compactly supported functions from $H^{1}({\mathbb{R}})$.
Note that $\mathfrak{t}_{q}(\varphi,\varphi)=\left(  L_{q}\varphi
,\varphi\right)  $ if $\varphi\in\mathcal{C}_{0}^{\infty}({\mathbb{R}})$, so
that if $L_{q}\geq0$, then both $\mathfrak{t}_{q}$ and $\mathfrak{t}_{q,I}$
are positive quadratic forms. Approximating $\varphi\in H^{1}_{comp}%
({\mathbb{R}})$ by functions from $\mathcal{C}_{0}^{\infty}(\mathbb{R})$, we
easily obtain that ${\mathfrak{t}}_{q}(\varphi,\varphi)\ge0$ for all
$\varphi\in H^{1}_{comp}({\mathbb{R}})$ as well.

It is easy to see that $\mathfrak{t}_{q,I}$ admits a closure, which has the
domain
\[
H_{0}^{1}(I)=\left\{  \psi\in H^{1}(I):\psi(a)=\psi(b)=0\right\}  .
\]
(See also Lemma 1.8 of \cite{SS:2003}.) It is a closed positive quadratic form
which will also be denoted $\mathfrak{t}_{q,I}$. Note that if $\varphi,\psi\in
H_{0}^{1}(I)$ and $\varphi_{0},\psi_{0}$ are their extensions on $\mathbb{R}$
by $0$, then $\varphi_{0},\psi_{0}\in H_{comp}^{1}({\mathbb{R}})$, and
\begin{equation}
\mathfrak{t}_{q,I}(\varphi,\psi)=\mathfrak{t}_{q}(\varphi_{0},\psi_{0}).
\label{E:tt0}%
\end{equation}

Let $L_{q,I}$ be the self-adjoint operator associated to $\mathfrak{t}_{q,I}$
by the Friedrichs construction. Clearly, $L_{q,I}$ has positive spectrum.
Moreover, it has compact resolvent, or, equivalently, discrete spectrum. (This
follows from the compactness of the imbedding of $H_{0}^{1}(I)$ into
$L^{2}(I)$; see also \cite{SS:2003}, where the asymptotics of the eigenvalues
is found.) By the min-max principle, the lowest eigenvalue of $L_{q,I}$ is
given by
\[
\lambda_{0}(I)=\inf\left\{  \mathfrak{t}_{q,I}(\varphi,\varphi):\varphi\in
H_{0}^{1}(I)\text{ and }\left\Vert \varphi\right\Vert _{L^{2}(I)}=1\right\}
\geq0
\]
and the infimum is achieved by a corresponding eigenfunction $h\in H_{0}%
^{1}(I)$.

\begin{lemma}
\label{lemma.dirichlet}Let $q$ be a real-valued distribution belonging to
$H_{\mathrm{loc}}^{-1}(\mathbb{R})$. If $L_{q}\geq0$, then $\lambda_{0}(I)>0$
for every bounded open interval $I$.
\end{lemma}

\begin{proof}
As $\lambda_{0}(I)\geq0$ it suffices to show that no bounded interval $I$ has
$\lambda_{0}(I)=0$. Suppose, on the contrary, that such an interval $I$
exists, and let $h\in H_{0}^{1}(I)$ be an $L^{2}$-normalized, real-valued
eigenfunction with the eigenvalue $0$, so, in particular, ${\mathfrak{t}%
}_{q,I}(h,h)=0$. Extend $h$ to a function $\eta=h_{0}\in H_{\mathrm{c}omp}%
^{1}(\mathbb{R})$ as above, i.e. by setting $\eta(x)=0$ if $x\in
\mathbb{R}\setminus I$. By \eqref{E:tt0},
we conclude that
$\mathfrak{t}_{q}\left(  \eta,\eta\right)  = \mathfrak{t}_{q,I}\left( h, h \right),$
hence $\mathfrak{t}_{q}\left(  \eta,\eta\right)  = 0$,
So, for any $\varphi\in\mathcal{C}_{0}^{\infty}(\mathbb{R})$ and
$t\in\mathbb{R}$,%
\[
\mathfrak{t}_{q}\left(  \eta+t\varphi,\eta+t\varphi\right)
=2t~\operatorname{Re}~\mathfrak{t}_{q}\left(  \eta,\varphi\right)
+t^{2}\text{ }\mathfrak{t}_{q}(\varphi,\varphi)
\]
is nonnegative due to positivity of ${\mathfrak{t}}_{q}$. It follows that
$~\operatorname{Re}~\mathfrak{t}_{q}(\eta,\varphi)=0$ for all $\varphi
\in\mathcal{C}_{0}^{\infty}(\mathbb{R})$, hence $\mathfrak{t}_{q}(\eta
,\varphi)=0$ also for all $\varphi\in\mathcal{C}_{0}^{\infty}(\mathbb{R})$. It
follows that $\eta$ solves the equation $L_{q}\eta=0$ and has a compact
suppport. By Lemma \ref{lemma.zero} we obtain $\eta(x)=0$ identically. This
contradicts the assumption that $\left\Vert h\right\Vert _{L^{2}(I)}=1$, and
the lemma is proved.
\end{proof}

\begin{corollary}
\label{corollary.dirichlet}Let $q$ be a real-valued distribution from
$H_{\mathrm{loc}}^{-1}(\mathbb{R})$. If $L_{q}\geq0$ and $y\in H_{\mathrm{loc}%
}^{1}(\mathbb{R})$ solves $L_{q}y=0$, then $y$ can have at most one zero.
\end{corollary}

\begin{proof}
If $y(a)=y(b)=0$ for $a<b$, then zero is a Dirichlet eigenvalue of $L_{q,I}$
with $I=(a,b)$, contradicting Lemma \ref{lemma.dirichlet}.\bigskip
\end{proof}

\begin{proposition}
\label{prop.pos}Let $q$ be a real-valued distribution which belongs to
$H_{\mathrm{loc}}^{-1}(\mathbb{R})$. If $L_{q}\geq0$, then $L_{q}y=0$ has a
strictly positive solution.
\end{proposition}

\begin{proof}
For $c\in\mathbb{R}\backslash\left\{  0\right\}  $, let $\left\{  y,u\right\}
\in H_{\mathrm{loc}}^{1}(\mathbb{R})\times W_{\mathrm{loc}}^{1,1}(\mathbb{R})$
be the unique solution of (\ref{eq.sys}) with $y(c)=0$ and $u(c)=1$. As $y$
satisfies $L_{q}y=0$, Corollary \ref{corollary.dirichlet} implies that
$y(0)\neq0$. Denote by $\left\{  y_{c},u_{c}\right\}  $ the scaled solution of
(\ref{eq.sys}),%
\begin{equation}
\label{E:yc}y_{c}(x)=y(x)/y(0),~~u_{c}(x)=u(x)/y(0).
\end{equation}
Then $y_{c}\in H_{\mathrm{loc}}^{1}(\mathbb{R})$ is the unique solution of
$L_{q}y=0$ with $y(0)=1$ and $y(c)=0$. Next, choose $c,c^{\prime}\in
\mathbb{R}$ so that $0<c^{\prime}<c$. Note that $w(x)=y_{c}(x)-y_{c^{\prime}%
}(x)$ is a solution of $L_{q}y=0$ with $w(0)=0$, but $w(c)>0$. By Corollary
\ref{corollary.dirichlet}, $w(x)>0$ for $x>0$ and by Lemma \ref{lemma.zero}
and Corollary \ref{corollary.dirichlet} $w(x)<0$ for $x<0$. It follows that,
on the half-line $c>0$, the map $c\mapsto y_{c}(x)$ is monotone decreasing for
any given $x<0$ and monotone increasing for any given $x>0$.

We wish to construct a positive solution of $L_{q}y=0$ by taking the limit
$c\rightarrow+\infty$. To this end denote by $(\tilde y_{1},\tilde u_{1})$ and
$(\tilde y_{2}, \tilde u_{2}) $ the fundamental solutions of (\ref{eq.sys}),
determined by $\tilde y_{1}(0)=1$, $\tilde u_{1}(0)=0$ and $\tilde
y_{2}(0)=0,\tilde u_{2}(0)=1$ respectively. By Lemma \ref{lemma.dirichlet},
$\tilde y_{2}(-1)\neq0.$ Hence, for any $\alpha\in\mathbb{R}$%
\[
z(x;\alpha):=\tilde y_{1}(x)+\frac{\alpha-\tilde y_{1}(-1)} {\tilde y_{2}%
(-1)}\tilde y_{2}(x)
\]
is the unique solution of $L_{q}y=0$ with $y(-1)=\alpha$ and $y(0)=1.$
Moreover, it follows that for any $x\in\mathbb{R},$ $\alpha\mapsto
z(x;\alpha)$ is continuous. (In fact, the map $\alpha\mapsto z(\cdot;\alpha)$
is an affine, hence continuous map ${\mathbb{R}}\to H^{1}_{\mathrm{l}%
oc}({\mathbb{R}})$.)

Note that $y_{c}(x)=z(x;y_{c}(-1))$. Now consider the solution $y_{n}$ of
$L_{q}y=0$, $y(0)=1$, $y(n)=0$. Then $\alpha_{n}=y_{n}(-1)$ is a strictly
positive, decreasing sequence. Let $\alpha_{\infty}=$ $\lim_{n\rightarrow
\infty}\alpha_{n}$. By the continuity of $z(x;\alpha)$ with respect to
$\alpha$ and the fact that $y_{n}(x)=z(x;\alpha_{n}),$ it then follows that
$z(x;\alpha_{\infty})=\lim_{n\rightarrow\infty}y_{n}(x)$ for any
$x\in\mathbb{R}$. We claim that $\alpha_{\infty}>0$ and that $z(x;\alpha
_{\infty})$ is positive. To prove this, first note that for any $n\geq1,$
$y_{n}(x)\geq0$ for all $x\leq n$ and hence $\lim_{n\rightarrow\infty}%
y_{n}(x)\geq0$ for all $x\in\mathbb{R}$. If $\alpha_{\infty}=0$, then
\[
0\leq\lim_{n\rightarrow\infty}y_{n}(-2)=\lim_{n\rightarrow\infty}%
z(-2;\alpha_{n})=z(-2;0)=y_{-1}(-2)<0,
\]
a contradiction. Hence $\alpha_{\infty}>0$. As a consequence, $z(x;\alpha
_{\infty})$ is a nonnegative solution of $L_{q}y=0$ and hence strictly
positive by Lemma \ref{lemma.zero}.
\end{proof}

\emph{Proof of Theorem \ref{thm.im1}}. In the statement of Theorem
\ref{thm.im1}, we have (ii)$\Rightarrow$(i) by Lemma \ref{lemma.straight}(a).
To show that (i)$\Rightarrow$(iii), we compute that, for $q=r^{\prime}+r^{2}$
and any $\varphi\in\mathcal{C}_{0}^{\infty}(\mathbb{R})$,
\begin{equation}
\left(  L_{q}\varphi,\varphi\right)  =\int\left\vert \varphi^{\prime}%
-r\varphi\right\vert ^{2}~dx\geq0. \label{eq.pos}%
\end{equation}
Finally, (iii)$\Rightarrow$(ii) by Proposition \ref{prop.pos}.%
\hfill$\Box$%

\section{The Image of the Miura Map}

In this section we prove Theorem \ref{thm.im2}. Recall that $B_{\beta}$
denotes the restriction of the Miura map to the Sobolev space $H^{\beta
}(\mathbb{R})$ for $\beta\geq0$.

We begin by considering the case $\beta=0$. In the light of Theorem
\ref{thm.im1}, we need to find necessary and sufficient conditions on a
potential $q\in H^{-1}(\mathbb{R})$ so that there exists a solution $r\in
L^{2}(\mathbb{R})$ of the Riccati equation $r^{\prime}+r^{2}=q$. Hartman
\cite{Hartman:1982}, chapter XI.7, Lemma 7.1 has studied this problem for
continuous $q$ and his arguments still apply in our more general setting.

\begin{lemma}
\label{lemma.Q}Suppose that $q\in H^{-1}(\mathbb{R})$ is a real-valued
distribution and that
\begin{equation}
\sup_{\left\vert T\right\vert >1}\left\vert \frac{1}{T}\int_{0}^{T}%
Q(x)~dx\right\vert <+\infty\label{eq.Q}%
\end{equation}
for an antiderivative $Q\in L_{\mathrm{loc}}^{2}(\mathbb{R})$ of $q$. Then
every solution $r\in L_{\mathrm{loc}}^{2}(\mathbb{R})$ of the Riccati equation
$r^{\prime}+r^{2}=q$ belongs to $L^{2}(\mathbb{R})$. Conversely, if $r\in
L^{2}(\mathbb{R})$, then every antiderivative $Q$ of $q=r^{\prime}+r^{2}$
satisfies (\ref{eq.Q}).
\end{lemma}

\begin{proof}
(i) Assume that $Q\in L_{\mathrm{loc}}^{2}(\mathbb{R})$, $Q^{\prime}=q$, and
$Q$ satisfies (\ref{eq.Q}). We need to show that for any solution $r\in
L_{\mathrm{loc}}^{2}(\mathbb{R})$ of $r^{\prime}+r^{2}=q$, the integrals
$\int_{0}^{\infty}r^{2}(s)~ds$ and $\int_{-\infty}^{0}r^{2}(s)~ds$ are finite.
Let us show that the first integral is finite. Since $Q$ is an antiderivative
of $q$, it follows that, for a constant $C$,%
\begin{equation}
r(x)+\int_{0}^{x}r^{2}(s)~ds=Q(x)+C. \label{eq.anti}%
\end{equation}
By assumption, the C\'{e}saro mean $T^{-1}\int_{0}^{T}\left(  ~\cdot~\right)
~dx$ of the right-hand side of (\ref{eq.anti}) is bounded as $T\rightarrow
+\infty$. We suppose that $\int_{0}^{x}r^{2}(s)~ds\rightarrow+\infty$ as
$x\rightarrow+\infty$ and obtain a contradiction as follows.

First note that if $f \in L_{\mathrm{loc}}^{1}(\mathbb{R})$ and $f(x)
\rightarrow+\infty$ as $x\rightarrow+\infty$ then the same holds for its
C\'{e}saro mean, i.e.
\begin{equation}
\frac{1}{T}\int_{0}^{T} f(x) ~dx \rightarrow+\infty\qquad{\mathrm{as}} \;\; T
\rightarrow+\infty. \label{C.means}%
\end{equation}

Hence, taking C\'{e}saro means of (\ref{eq.anti}) we see that, if $\int
_{0}^{x}r^{2}(s)~ds \rightarrow+\infty$ as $x\rightarrow+\infty$, then
$T^{-1}\int_{0}^{T}~r(x)~dx\rightarrow-\infty$ as $x\rightarrow+\infty$.
Moreover, there is a $T_{0}$ so that for all $T>T_{0}$,
\[
-\frac{2}{T}\int_{0}^{T}~r(x)~dx\geq\frac{1}{T}\int_{0}^{T}~\left(  \int
_{0}^{x}~r^{2}(s)~ds\right)  ~dx>0.
\]
By the Cauchy-Schwarz inequality,
\[
-T^{-1}\int_{0}^{T}~r(x)~dx\leq T^{-1/2}\left(  \int_{0}^{T}~r^{2}%
(x)~dx\right)  ^{1/2}%
\]
so that%
\[
4T\int_{0}^{T}~r^{2}(x)~dx\geq\left[  \int_{0}^{T}\left(  \int_{0}^{x}%
~r^{2}(s)~ds\right)  dx\right]  ^{2}.
\]
Setting $I(T)=\int_{0}^{T}\int_{0}^{x}~r^{2}(s)~ds~dx$, we have that%
\[
4T~I^{\prime}(T)\geq I(T)^{2}%
\]
from which it follows by integration that%
\[
\frac{1}{I(T_{0})}-\frac{1}{I(T)}\geq\frac{1}{4}\log\left(  T/T_{0}\right)  .
\]
This contradicts that, by (\ref{C.means}), $I(T)\rightarrow+\infty$ as
$T\rightarrow+\infty$. A similar argument shows that $\int_{-\infty}^{0}%
r^{2}(s)~ds$ is finite. Hence $r\in L^{2}(\mathbb{R})$.

(ii) If, on the other hand, $q=r^{\prime}+r^{2}$ for $r\in L^{2}(\mathbb{R})$,
then the function $Q(x)=r(x)+\int_{0}^{x}r^{2}(s)~ds$ satisfies (\ref{eq.Q})
by the Cauchy-Schwarz inequality applied in $\left[  0,T\right]  $.
\end{proof}

\begin{corollary}
\label{corollary.Q}Suppose that $q\in\operatorname{Im}(B_{0})$, i.e.,
$q=r^{\prime}+r^{2}$ for some $r\in L^{2}(\mathbb{R})$. If $u\in
L_{\mathrm{loc}}^{2}(\mathbb{R})$ solves the Riccati equation $u^{\prime
}+u^{2}=q$, then $u\in L^{2}(\mathbb{R})$.
\end{corollary}

\begin{proposition}
\label{prop.image2}A real-valued distribution $q\in H^{-1}(\mathbb{R})$
belongs to $\operatorname{Im}(B_{0})$ if and only if\newline(i) $L_{q}\geq0$,
and\newline(ii) $q$ can be presented as $q=f^{\prime}+g$ for real-valued
functions $f\in L^{2}(\mathbb{R})$ and $g\in L^{1}(\mathbb{R})$.
\end{proposition}

\begin{proof}
Suppose that $q\in\operatorname{Im}(B_{0})$, i.e., $q=r^{\prime}+r^{2}$ for
some $r\in L^{2}(\mathbb{R})$. Then $L_{q}\geq0$ by Theorem \ref{thm.im1} and
$q=f^{\prime}+g$ with $f=r$, $g=r^{2}$. On the other hand, suppose that $q\in
H^{-1}(\mathbb{R})$ with $L_{q}\geq0$, and $q=f^{\prime}+g$ for $f\in
L^{2}(\mathbb{R})$ and $g\in L^{1}(\mathbb{R})$. By Theorem \ref{thm.im1},
$q\in\operatorname{Im}(B)$, so $q=r^{\prime}+r^{2}$ for some $r\in
L_{\mathrm{loc}}^{2}(\mathbb{R})$. The antiderivative $Q(x)=f(x)+\int_{0}%
^{x}g(s)~ds$ obeys the condition (\ref{eq.Q}), so $r\in L^{2}(\mathbb{R})$ by
Lemma \ref{lemma.Q}.
\end{proof}

\begin{proposition}
\label{prop.sicktop}The set $\operatorname{Im}(B_{0})$ has no interior points,
and hence is not open in $H^{-1}(\mathbb{R})$. Further, the set
$\operatorname{Im}(B_{0})$ is not closed in $H^{-1}(\mathbb{R})$.
\end{proposition}

\begin{proof}
First we show that $\operatorname{Im}(B_{0})$ has empty interior. If
$q\in\mathcal{C}_{0}^{\infty}(\mathbb{R})\cap\operatorname{Im}(B_{0})$, we can
perturb $q$ by a small potential well far separated from the support of $q$
and create a bound state. More precisely, for $\varepsilon>0$, let
\[
v_{\varepsilon}(x)=\left\{
\begin{array}
[c]{cll}%
-\varepsilon &  & \left\vert x\right\vert <1/(2\varepsilon)\\
0 &  & \left\vert x\right\vert \geq1/(2\varepsilon)
\end{array}
\right.  .
\]
Observe that $\int v_{\varepsilon}(x)dx=-1$ but $\left\Vert v_{\epsilon
}\right\Vert _{L^{2}(\mathbb{R})}=\varepsilon$. Suppose that $q\in
\mathcal{C}_{0}^{\infty}(\mathbb{R})\cap\operatorname{Im}(B_{0})$ with support
contained in $\left[  -a,a\right]  $. Let $w_{\varepsilon}(x)=v_{\varepsilon
}(x-2a-2\varepsilon^{-1})$; then $w_{\varepsilon}$ has support disjoint from
the one of $q$. By choosing $\varepsilon$ sufficiently small, we can assure
that the potential $q_{\varepsilon}=q+w_{\varepsilon}$ is close to $q$ in
$L^{2}(\mathbb{R})$ norm. \ Let $\chi\in\mathcal{C}_{0}^{\infty}(\mathbb{R})$
be a nonnegative function with $\chi(x)=1$ for $\left\vert x\right\vert
<1/(2\varepsilon)$, $\chi(x)=0$ for $\left\vert x\right\vert >\varepsilon
^{-1}$, and $\left\vert \chi^{\prime}(x)\right\vert \leq3\varepsilon$.
Finally, let $\eta(x)=$ $\chi\left(  x-2a-2\varepsilon^{-1}\right)  $. Then%
\begin{align*}
\left(  L_{q_{\varepsilon}}\eta,\eta\right)   &  =\int\left\vert \eta^{\prime
}(x)\right\vert ^{2}-1 \leq18\varepsilon-1.
\end{align*}
Hence, by Theorem \ref{thm.im1}, $q_{\varepsilon}\notin\operatorname{Im}(B)$
for $0 < \varepsilon< 1/18$. Since $\mathcal{C}_{0}^{\infty}(\mathbb{R})$
is norm-dense in $H^{-1}(\mathbb{R})$, this shows that $\operatorname{Im}%
(B_{0})$ contains no open neighborhood in the norm topology of $H^{-1}%
(\mathbb{R})$.

Next, we show that $\operatorname{Im}(B_{0})$ is not closed. Suppose that $q$
is any nonnegative function with $q\in L^{2}(\mathbb{R})$ but $q\notin
L^{1}(\mathbb{R})$. We may approximate $q$ by nonnegative potentials $q_{k}%
\in\mathcal{C}_{0}^{\infty}(\mathbb{R})$ so that $q_{k}\rightarrow q$ in
$L^{2}(\mathbb{R})$. Hence, $q_{k}\rightarrow q$ in $H^{-1}(\mathbb{R})$ and
by Proposition \ref{prop.image2}, $q_{k}\in\operatorname{Im}(B_{0})$ for any
$k\ge1.$ Moreover, since $L_{q}\geq0$, it follows from Theorem \ref{thm.im1}
that $q\in\operatorname{Im}(B)$. On the other hand, as $q \ge0,$ no
antiderivative $Q$ of $q$ satisfies condition (\ref{eq.Q}). Thus
$\operatorname{Im}(B_{0})$ is not closed in the norm topology on
$H^{-1}(\mathbb{R})$.
\end{proof}

The image of $B_{0}$ can also be characterized by a `special integral' of $q$.
Let $\left\{  \chi_{n}\right\}  _{n\geq1}$ be a sequence of nonnegative
$\mathcal{C}_{0}^{\infty}(\mathbb{R})$ functions with (i) $\chi_{n}(x)=1$ for
$\left\vert x\right\vert \leq n$, (ii) $\chi_{n}(x)=0$ for $\left\vert
x\right\vert \geq n+1$, and (iii) $\left\vert \chi_{n}^{\prime}(x)\right\vert
\leq2$ for all $x\in\mathbb{R}$. Given $q\in H^{-1}(\mathbb{R})$, we define
the special integral of $q$, denoted $\left[  q\right]  $, to be the number
$\lim_{n\rightarrow\infty}(q,\chi_{n})$ if this limit exists and is finite.
One easily checks that $\left[  q\right]  $ is well-defined, i.e., does not
depend on the choice of sequence $\left\{  \chi_{n}\right\}  _{n\geq1}$
satisfying properties (i), (ii), and (iii) above. If $q\in\operatorname{Im}%
(B_{0})$ then, for any $r\in B_{0}^{-1}(q)$,%
\begin{align}
\lim_{n\rightarrow\infty} \left(  q,\chi_{n}\right)   &  =\lim_{n\rightarrow
\infty}\left\{  \left(  -r,\chi_{n}^{\prime}\right)  +\left(  r^{2},\chi
_{n}\right)  \right\}  =\left\Vert r\right\Vert _{L^{2}(\mathbb{R})}^{2}
\label{eq.r}%
\end{align}
which shows that $\left[  q\right]  \geq0$ on $\operatorname{Im}(B_{0})$ with
$\left[  q\right]  =0$ if and only if $q=0$. Moreover, if $r_{1}$ and $r_{2}$
belong to $B^{-1}(q)$, then $\left\Vert r_{1}\right\Vert _{L^{2}(\mathbb{R}%
)}=\left\Vert r_{2}\right\Vert _{L^{2}(\mathbb{R})}$.

The special integral has the following properties:

\begin{enumerate}
\item[(a)] $\left[  f\right]  =\int_{\mathbb{R}}f~dx$ if $f\in L^{1}%
(\mathbb{R})\subset H^{-1}(\mathbb{R})$;

\item[(b)] $\operatorname*{Dom}\left(  \left[  ~\cdot~\right]  \right)  $ is a
linear subspace in $H^{-1}(\mathbb{R})$, and $f\mapsto\left[  f\right]  $ is linear;

\item[(c)] $\left[  f^{\prime}\right]  =0$ for any $f\in L^{2}(\mathbb{R})$;

\item[(d)] If $f \in L^{2}(\mathbb{R})$, then $\left[  f \right]  $ exists if
and only if $f$ is conditionally integrable, i.e. the limit $\lim
_{T\rightarrow\infty} \int_{-T}^{T} f(x) ~dx$ exists. In this case, the limit
equals $\left[  f \right]  .$
\end{enumerate}

Using the special integral we can give an alternative characterization of
$\operatorname{Im}(B_{0})$.

\begin{theorem}
\label{thm.special}A real-valued distribution $q\in H^{-1}(\mathbb{R})$
belongs to $\operatorname{Im}(B_{0})$ if and only if:\newline(i) $L_{q}\geq0$,
and\newline(ii) $\left[  q\right]  $ exists. \newline Moreover, for any $q
\in\operatorname{Im}(B_{0}), $ one has $\left[  q\right]  \ge0. $
\end{theorem}

\begin{proof}
First, suppose that $L_{q}\geq0$ and $\left[  q\right]  $ exists. To prove
that $q\in\operatorname{Im}(B_{0})$, it suffices by Lemma \ref{lemma.Q} to
show that $q$ has an antiderivative $Q$ with bounded C\'{e}saro means. By Lemma
\ref{lemma.rep}, any $q\in H^{-1}(\mathbb{R})$ may be written $q=f^{\prime}+g$
for $f$ and $g$ belonging to $L^{2}(\mathbb{R})$. We can therefore take
$Q(x)=f(x)+G(x)$ where $G(x)=\int_{0}^{x}g(s)~ds$. We will use condition (ii)
on $q$ to show that $G$ is bounded. Since $\left[  f^{\prime}\right]  =0$ for
any $f\in L^{2}(\mathbb{R})$, we have $\left[  q\right]  =\left[  g\right]  $
and $\left[  g\right]  $ exists. Since, also, $g\in L^{2}(\mathbb{R})$, the
existence of $\left[  g\right]  $ implies that $g$ is conditionally
integrable. Thus, $\lim_{n\rightarrow\infty}\alpha_{n}$ exists where
$\alpha_{n}:=\int_{-n}^{n}g(x)~dx$. This is equivalent to the existence of
$\lim_{a\rightarrow+\infty}\int_{-a}^{a}g(x)~dx$ if $g\in L^{2}(\mathbb{R})$.
We need to show that the numbers $\alpha_{n}^{+}=\int_{0}^{n}~g(x)~dx$ and
$\alpha_{n}^{-}=\int_{-n}^{0}~g(x)~dx$ are also bounded. Let $\left\{
\eta_{n}\right\}  _{n\geq1}$ be a sequence of $\mathcal{C}_{0}^{\infty
}(\mathbb{R})$ functions with $0\leq\eta_{n}(x)\leq1$, $\eta_{n}(x)=1$ for
$x\in\left[  0,n\right]  $, $\eta_{n}(x)=0$ for $x\in\mathbb{R}\backslash
\lbrack-1,n+1]$, and $\left\vert \eta_{n}^{\prime}(x)\right\vert \leq2$. Since
$(L_{q}\eta_{n},\eta_{n})\geq0$,
\[
-\left(  f, 2 \eta_{n}\eta_{n}^{\prime}\right)  +\int g\eta_{n}^{2} ~ dx
\geq-\left\Vert \eta_{n}^{\prime}\right\Vert _{L^{2}(\mathbb{R})}^{2}.
\]
Since $\left\Vert \eta_{n}^{\prime}\right\Vert _{L^{2}(\mathbb{R})} \le4$ and
\[
\left\vert \left(  f, 2 \eta_{n}\eta_{n}^{\prime} \right)  \right\vert \le4
\int_{-1}^{0} \left\vert f(x)\right\vert ~ dx + 4 \int_{n}^{n+1} \left\vert
f(x)\right\vert ~ dx \le8 \left\Vert f \right\Vert _{L^{2}(\mathbb{R})}
\]
as well as
\[
\int g\eta_{n}^{2} ~ dx = \alpha_{n}^{+} + \int_{-1}^{0} g(x) \eta_{n}^{2}(x)
~ dx + \int_{n}^{n+1} g(x) \eta_{n}^{2}(x) ~ dx \le\alpha_{n}^{+} + 2
\left\Vert g \right\Vert _{L^{2}(\mathbb{R})}
\]
it then follows that $\alpha_{n}^{+}\geq-C$ with $C$ independent of $n$. A
similar argument shows that $\alpha_{n}^{-}\geq-C$ with $C$ independent of
$n$. As $\alpha_{n}=\alpha_{n}^{+}+\alpha_{n}^{-}$ we conclude that the
sequences $\left\{  \alpha_{n}^{+}\right\}  $ and $\left\{  \alpha_{n}%
^{-}\right\}  $ are both bounded, so $G$ is bounded. Since $L_{q}\geq0$ we
have $r^{\prime}+r^{2}=q$ for some $r\in L_{\mathrm{loc}}^{2}(\mathbb{R})$ by
Theorem \ref{thm.im1}, and applying Lemma \ref{lemma.Q} we conclude that $r\in
L^{2}(\mathbb{R})$. Hence $q\in\operatorname{Im}(B_{0})$.

On the other hand, if $q\in\operatorname{Im}(B_{0})$, then $L_{q}\geq0$ by
Theorem \ref{thm.im1} and $q = r^{\prime} + r^{2}$ for some $r \in
L^{2}(\mathbb{R}),$ hence by properties (c) and (d) of the special integral,
$\left[  q\right]  $ exists and $\left[  q\right]  = \left\Vert r \right\Vert
_{L^{2}(\mathbb{R})}^{2} \ge0$.
\end{proof}

\begin{corollary}
\label{corollary.odd} An odd distribution $q \in H^{-1}(\mathbb{R})$ cannot be
in $\operatorname{Im}(B_{0})$ unless $q \equiv0.$
\end{corollary}

\begin{remark}
Generally, the condition $L_{q}\geq0$ can be considered as a weak form of
positivity for $q$. If it is satisfied then the existence of the special
integral $[q]$ for $q\in H^{-1}{(\mathbb{R}})$ implies much stronger
existence-of-limit type results. For example, let us take any family of
functions $\chi_{T_{1},T_{2}}\in\mathcal{C}_{0}^{\infty}({\mathbb{R}})$,
$T_{1},T_{2}\in{\mathbb{R}}$, such that $\chi_{T_{1},T_{2}}=1$ on
$[-T_{1},T_{2}]$, $\chi_{T_{1},T_{2}}=0$ on ${\mathbb{R}}\setminus
(-T_{1}-1,T_{2}+1)$, and the derivatives $\chi_{T_{1},T_{2}}^{\prime}(x)$ are
uniformly bounded. Then $\L _{q}\geq0$ and existence of $[q]$ imply the
existence of the limit
\[
\lim_{T_{1},T_{2}\rightarrow+\infty}\langle q,\chi_{T_{1},T_{2}}\rangle,
\]
which in case $q\in L^{2}({\mathbb{R}})$ is equivalent to the existence of the
limit
\[
\lim_{T_{1},T_{2}\rightarrow+\infty}\int_{-T_{1}}^{T_{2}}q(x)dx.
\]
To prove the above statements we can, for example, use Theorem
\ref{thm.special} to find \newline$r\in L^{2}({\mathbb{R}})$, such that
$q=r^{\prime}+r^{2}$, and the result easily follows.
\end{remark}

\medskip We now consider the restriction $B_{\beta}:H^{\beta}(\mathbb{R}%
)\rightarrow H^{\beta-1}(\mathbb{R})$ for $\beta>0$.

\begin{lemma}
\label{lemma.regular}Let $\beta\geq0$. If $q\in\operatorname{Im}(B_{0})\cap
H^{\beta-1}(\mathbb{R})$ and $r\in L_{\mathrm{loc}}^{2}(\mathbb{R})$ is a
solution of the Riccati equation $r^{\prime}+r^{2}=q$, then $r\in H^{\beta
}(\mathbb{R})$.
\end{lemma}

\begin{proof}
By Corollary \ref{corollary.Q}, the result holds for $\beta=0$. Hence, it
suffices to prove that in case the claimed result holds for a given $\beta
_{0}\geq0$, it also holds for any $\beta\in\left[  \beta_{0},\beta
_{0}+1/4\right]  $. So, assume that $q := r^{\prime}+r^{2}\in H^{\beta
-1}(\mathbb{R})$ with $\beta_{0} < \beta\le\beta_{0} +\frac{1}{4}$ and $r\in
H^{\beta_{0}}(\mathbb{R})$. Then $r^{2}$ belongs to $H^{-1/2-\delta
}(\mathbb{R})$, $H^{2\beta_{0}-1/2}(\mathbb{R})$, $H^{1/2-\delta}(\mathbb{R}%
)$, or $H^{\beta_{0}}(\mathbb{R})$ respectively when $\beta_{0}=0$,
$0<\beta_{0}<1/2$, $\beta_{0}=1/2$, or $\beta_{0}>1/2$ (see (\ref{eq.H1}%
)-(\ref{eq.H4})). In the first and last cases, $\delta>0$ can be chosen
arbitrarily small. Then $r^{\prime}=q-r^{2}$ is in $H^{s-1}(\mathbb{R})$ with
$s=\min(\beta,1/2-\delta),~\min(\beta,2\beta_{0}+1/2)$, $\min(\beta
,3/2-\delta)$ or $\min(\beta,\beta_{0}+1)$ respectively. As $r \in
L^{2}(\mathbb{R})$ it then follows that $r \in H^{s}(\mathbb{R})$ and since
$\beta_{0}\leq\beta\leq\beta_{0}+1/4$ implies $s=\beta$ in all cases, we get
the desired result.
\end{proof}

\emph{Proof of Theorem \ref{thm.im2}.} First, suppose that $q\in H^{\beta
-1}(\mathbb{R})$ satisfies conditions (i) and (ii) of Theorem \ref{thm.im2}.
From the trivial inclusion $H^{\beta}(\mathbb{R})\subset H^{0}(\mathbb{R})$
and Proposition \ref{prop.image2}, it follows that $q\in\operatorname{Im}%
(B_{0})$, i.e., $q=r^{\prime}+r^{2}$ for some function $r\in L^{2}%
(\mathbb{R})$. Applying Lemma \ref{lemma.regular} we see that $r\in H^{\beta
}(\mathbb{R})$, so $q\in\operatorname{Im}(B_{\beta})$ as claimed. Second, if
$q\in\operatorname{Im}(B_{\beta})$, then $L_{q}\geq0$ by Theorem \ref{thm.im1}
and $q=f^{\prime}+g$ with $f=r\in L^{2}(\mathbb{R})$ and $g=r^{2}\in
L^{1}(\mathbb{R})$.%
\hfill$\Box$%

\section{Geometry of the Miura Map}

In this section, we prove Theorem \ref{thm.im3}. According to Proposition
\ref{prop.homeo}, $B^{-1}(q)$ is homeomorphic to the set of positive solutions
$y$ of $L_{q}y=0$ with $y\in H_{\mathrm{loc}}^{1}(\mathbb{R})$ and $y(0)=1$.
As before we denote this set by $\operatorname*{Pos}(q)$. We will show that
$\operatorname*{Pos}(q)$ is either a point or homeomorphic to a line segment.

Suppose that $\operatorname*{Pos}(q)$ is nonempty and choose $y_{1}%
\in\operatorname*{Pos}(q)$. Using the Wronskian we can find another solution%
\[
y_{2}(x)=y_{1}(x)\int_{0}^{x}y_{1}(s)^{-2}~ds.
\]
The general solution to $L_{q}y=0$ is then%
\[
y(x)=y_{1}(x)\left(  c_{1}+c_{2}F(x)\right)
\]
where%
\[
F(x)=\int_{0}^{x}y_{1}(s)^{-2}ds.
\]
Observe that $F$ is a monotone increasing function with $F(0)=0$. If we define
numbers $m_{\pm}\in(0,+\infty]$ by
\begin{equation}
m_{+}=\lim_{x\rightarrow+\infty}F(x) \label{eq.m+}%
\end{equation}
and
\begin{equation}
m_{-}=-\lim_{x\rightarrow-\infty}F(x), \label{eq.m-}%
\end{equation}
then $F$ takes values in $\left(  -m_{-},m_{+}\right)  $. We will set
$m_{+}^{-1}=0$ if $m_{+}=+\infty$, and similarly for $m_{-}^{-1}$. The
conditions $y(0)=1$ and $y(x)>0$ for all $x$ determine that any $y\in
\operatorname*{Pos}(q)$ is written%
\[
y(x)=y_{1}(x)\left(  1+cF(x)\right)
\]
with $c\in\left[  -m_{+}^{-1},m_{-}^{-1}\right]  $. Letting%
\begin{align}
y_{+}(x)  &  =y_{1}(x)\left(  1-m_{+}^{-1}F(x)\right) \label{eq.y+}\\
y_{-}(x)  &  =y_{1}(x)\left(  1+m_{-}^{-1}F(x)\right)  \label{eq.y-}%
\end{align}
we see that%
\begin{align}
y_{+}(x)  &  \leq y(x)\leq y_{-}(x),~~x>0\label{eq.bd1}\\
y_{-}(x)  &  \leq y(x)\leq y_{+}(x),~~x<0 \label{eq.bd2}%
\end{align}
for any $y\in\operatorname*{Pos}(q)$, and%
\[
\operatorname*{Pos}(q)=\left\{  \theta y_{+}+(1-\theta)y_{-}:\theta\in\left[
0,1\right]  \right\}  .
\]
Thus, either (i) $m_{+}=m_{-}=+\infty$, $y_{+}=y_{-}$ and $\operatorname*{Pos}%
(q)$ consists of a single element, or (ii) at least one of $m_{\pm}$ is
finite, $y_{+}\neq y_{-}.$ Noting that $\theta\mapsto\theta y_{+} + (1-\theta)
y_{-}$ is a continuous map from $[0,1]$ to the Hausdorff space
$H_{\mathrm{loc}}^{1}(\mathbb{R})$ we see that $\operatorname*{Pos}(q)$ is
homeomorphic to the interval $\left[  0,1\right]  $. We have proved:

\begin{lemma}
\label{lemma.dichotomy}Suppose that $q\in H_{\mathrm{loc}}^{-1}(\mathbb{R})$
is a real-valued distribution and $L_{q}\geq0$. Then $\operatorname*{Pos}(q)$
is either a point or homeomorphic to a line segment.
\end{lemma}

Next, we show that the sets $E_{1}$ and $E_{2}$ defined in (\ref{eq.e1}) and
(\ref{eq.e2}) are both dense in $H_{\mathrm{loc}}^{-1}(\mathbb{R})$. We begin
with a simple lemma which will be useful in the proof of Theorem \ref{thm.im3}.

\begin{lemma}
\label{lemma.e2}There exists a family of potentials $\left\{  w_{\varepsilon
}\right\}  _{\varepsilon>0}$ contained in $\mathcal{C}_{0}^{\infty}%
(\mathbb{R})\cap E_{2}$ so that (i) $\operatorname*{supp}(w_{\varepsilon
})\subset\left[  -\varepsilon^{-1},\varepsilon^{-1}\right]  $ and (ii)
$\left\Vert w_{\varepsilon}\right\Vert _{H^{\beta}(\mathbb{R})}\rightarrow0$
as $\varepsilon\downarrow0$ for any $\beta\in\mathbb{R}$.
\end{lemma}

\begin{proof}
Let $y\in\mathcal{C}^{\infty}(\mathbb{R})$ with $y(x)=1$ for $x<-1$, $y(x)=x$
for $x>1$, and $y(x)>0$ for any $x \in\mathbb{R}$. The potential
$w(x)=y^{\prime\prime}(x)/y(x)$ has $y$ as a positive solution of $L_{w}y=0$
and $w=B(r)$ with $r=y^{\prime}/y$. Since $\int_{0}^{\infty}y(s)^{-2}%
ds<\infty$, it follows from the remarks preceding Lemma \ref{lemma.dichotomy}
that $w\in E_{2}$. Now let $y_{\varepsilon}(x)=y\left(  \varepsilon x\right)
$ and $w_{\varepsilon}(x)=y_{\varepsilon}^{\prime\prime}(x)/y_{\varepsilon
}(x)$. Then $w_{\varepsilon}\in E_{2}$ with support in $\left[  -\varepsilon
^{-1},\varepsilon^{-1}\right]  $, proving (i). To prove (ii), note that
$w_{\varepsilon}(x)=\varepsilon^{2}w\left(  \varepsilon x\right)  $ so that
for any nonnegative integer $j$,
\[
\left\Vert \partial_{x}^{j}w_{\varepsilon}\right\Vert _{L^{2}(\mathbb{R})}%
^{2}=\varepsilon^{3+2j}\left\Vert \partial_{x}^{j}w\right\Vert _{L^{2}%
(\mathbb{R})}^{2}.
\]
Since $\left\Vert u\right\Vert _{H^{\alpha}(\mathbb{R})}\leq\left\Vert
u\right\Vert _{H^{\beta}(\mathbb{R})}$ for $\alpha<\beta$ and $u\in H^{\beta
}(\mathbb{R})$, this shows that $\left\Vert w_{\varepsilon}\right\Vert
_{H^{\beta}(\mathbb{R})}\rightarrow0$ as $\varepsilon\downarrow0,$ for any
$\beta\in\mathbb{R}$.
\end{proof}

\begin{lemma}
\label{lemma.dense1}$B(\mathcal{C}_{0}^{\infty}(\mathbb{R}))$ is dense in
$\operatorname{Im}(B)$ and $B(\mathcal{C}_{0}^{\infty}(\mathbb{R}))\subset
E_{1}$.

\begin{proof}
Let $q\in\operatorname{Im}(B)$ and let $r\in B^{-1}(q)$. Let $\left\{
r_{n}\right\}  \in\mathcal{C}_{0}^{\infty}(\mathbb{R})$ with $r_{n}\rightarrow
r$ in $L_{\mathrm{loc}}^{2}(\mathbb{R})$. By the continuity of the Miura map,
$B(r_{n})\rightarrow B(r)$ in $H_{\mathrm{loc}}^{-1}(\mathbb{R})$. Thus
$B\left(  \mathcal{C}_{0}^{\infty}(\mathbb{R})\right)  $ is dense in
$\operatorname{Im}(B)$. If $q=B(r)$ for $r\in\mathcal{C}_{0}^{\infty
}(\mathbb{R})$, then $\operatorname*{Pos}(q)$ contains the element
$y_{1}(x)=\exp\left(  \int_{0}^{x}r(s)~ds\right)  $ which is bounded above and
below by strictly positive constants. It follows that $m_{+}=m_{-}=+\infty$
(see (\ref{eq.m+}) and (\ref{eq.m-})), By the analysis of positive solutions
preceding Lemma \ref{lemma.dichotomy}, $y_{+}=$ $y_{-}$ and
$\operatorname*{Pos}(q)$ consists of a single point. Hence $B(\mathcal{C}%
_{0}^{\infty}(\mathbb{R}))\subset E_{1}$.
\end{proof}
\end{lemma}

On the other hand:

\begin{lemma}
\label{lemma.dense2}$E_{2}$ is dense in $\operatorname{Im}(B)$.
\end{lemma}

\begin{proof}
Since $B(\mathcal{C}_{0}^{\infty}(\mathbb{R}))$ is dense in $\operatorname{Im}%
(B)$, it suffices to show that for any $q\in B(\mathcal{C}_{0}^{\infty
}(\mathbb{R}))$ there is a sequence of elements $q_{n}$ from $E_{2}$ with
$q_{n}\rightarrow q$ in $H_{\mathrm{loc}}^{-1}(\mathbb{R})$ as $n\rightarrow
\infty$. Suppose that $q\in B(\mathcal{C}_{0}^{\infty}(\mathbb{R}))$ with
support in $\left[  -a,a\right]  $ for $a>0$ and consider the sequence%
\[
q_{n}=q+v_{n}%
\]
for $n\geq1$, where
\[
v_{n}(x)=w_{1/n}(x-a-2n)
\]
and $w_{\varepsilon}$ is the family constructed in Lemma \ref{lemma.e2}. Then
$v_{n} \rightarrow0$ for $n \rightarrow\infty$ in $H^{\beta}(\mathbb{R})$ for
any $\beta\in\mathbb{R}.$ Let $p \in H_{\mathrm{loc}}^{1}(\mathbb{R})$ be the
unique positive solution to $L_{q}p=0$ with $p(0)=1$; note that $p(x)$ is
constant away from the support of $q$. If $y_{\varepsilon}$ is the positive
solution for $w_{\varepsilon}$ constructed in Lemma \ref{lemma.e2}, it is
easily seen that the function%
\begin{equation}
z_{n}(x)=\left\{
\begin{array}
[c]{lc}%
p(x), & x<a+1,\\
p(a+1)y_{1/n}\left(  x-a-2n\right)  , & x\geq a+1,
\end{array}
\right.  \label{eq.yn}%
\end{equation}
is a positive solution to $L_{q_{n}}y=0$. It follows from (\ref{eq.yn}) and
the fact that $y_{1/n}(x)=x/n$ for $x$ large and positive that $\int
_{0}^{\infty}z_{n}(s)^{-2}ds<\infty$. Thus, by the analysis of positive
solutions preceding Lemma \ref{lemma.dichotomy}, $q_{n}\in E_{2}$. Since
$v_{n}\rightarrow0$ in $H^{\beta}(\mathbb{R})$ for any $\beta\in\mathbb{R}$,
$q_{n}-q\rightarrow0$ in $H_{\mathrm{loc}}^{-1}(\mathbb{R})$.
\end{proof}

\emph{Proof of Theorem \ref{thm.im3}.} That $\operatorname{Im}(B)=E_{1}\cup
E_{2}$ follows from Lemma \ref{lemma.dichotomy} and Proposition
\ref{prop.homeo}. The density statements were proved in Lemmas
\ref{lemma.dense1} and \ref{lemma.dense2}.
\hfill$\Box$%
\newline

We close this section with some further remarks on the dichotomy of the Miura
map. First, we give a version of Theorem \ref{thm.im3} for the restriction of
the Miura map to $H^{\beta}(\mathbb{R})$, $\beta\geq0$.

\begin{theorem}
$\operatorname{Im}(B_{\beta})=E_{1,\beta}\cup E_{2,\beta}$ where $E_{j,\beta
}=E_{j}\cap\operatorname{Im}(B_{0})\cap H^{\beta-1}(\mathbb{R}),$ $j = 1,2$.
Moreover $E_{j,\beta}$ is dense in $\operatorname{Im}(B_{\beta})$ for $j=1,2$.
\end{theorem}

\begin{proof}
The first statement follows from the fact, established in Theorem
\ref{thm.im2}, that $\operatorname{Im}(B_{\beta})=\operatorname{Im}(B_{0})\cap
H^{\beta-1}(\mathbb{R})$. The proofs of Lemmas \ref{lemma.dense1} and
\ref{lemma.dense2} can be adapted with trivial changes to show the density of
$E_{j,\beta}$ in $\operatorname{Im}(B_{\beta})$.
\end{proof}

Finally, let
\begin{equation}
\lambda_{0}(q)=\inf\left\{  \left(  L_{q}\varphi,\varphi\right)  :\varphi
\in\mathcal{C}_{0}^{\infty}(\mathbb{R}),~\left\Vert \varphi\right\Vert
_{L^{2}(\mathbb{R})}=1\right\}  , \label{eq.lambda}%
\end{equation}
or, equivalently,
\begin{equation}
\lambda_{0}(q)=\inf\left\{  \frac{\mathfrak{t}_{q}(\psi,\psi)}{(\psi,\psi)}:
\psi\in H^{1}_{\mathrm{comp}}(\mathbb{R})\setminus\{0\}\right\}  ,
\label{eq.qlambda}%
\end{equation}
and define the sets%
\[
E_{\bullet}=\left\{  q\in\operatorname{Im}(B):\lambda_{0}(q)=0\right\}
\]
and
\[
E_{>}=\left\{  q\in\operatorname{Im}(B):\lambda_{0}(q)>0\right\}  .
\]
If $q$ has compact support, it is clear that $\lambda_{0}(q)=0$ since we can
choose test functions whose support is disjoint from the support of $q$ and
\[
\inf\left\{  \left\Vert \varphi^{\prime}\right\Vert ^{2}:\varphi\in
\mathcal{C}_{0}^{\infty}(\mathbb{R}),~\left\Vert \varphi\right\Vert
_{L^{2}(\mathbb{R})}=1\right\}  =0.
\]
Note that the map $E_{\bullet}\times\mathbb{R}_{>0}\rightarrow E_{>}$ given by%
\[
\left(  q,c\right)  \mapsto q+c
\]
is a continuous, bijective map onto $E_{>}$.

\begin{theorem}
\label{thm.nofold}(i) $E_{>}\subset E_{2}$ and $E_{1}\subseteq E_{\bullet}%
$.\newline(ii) $E_{2}\cap E_{\bullet}\neq\emptyset$, i.e., $E_{\bullet}$ is
not a fold of the dichotomy. Moreover, $E_{\bullet}$ is dense in
$\operatorname{Im}(B)$, and $E_{1}$ and $E_{2}\cap E_{\bullet}$ are dense in
$E_{\bullet}$.
\end{theorem}

\begin{remark}
The fact that $E_{2}\cap E_{\bullet}\neq\emptyset$ has already been observed
by Murata \cite{Murata:1986}, Remark to Theorem 2.2, in his investigation of
critical and subcritical potentials -- see Appendix \ref{app.work}.
\end{remark}

The proof of Theorem \ref{thm.nofold} will rely on the following proposition.

\begin{proposition}
\label{prop.2pos} Suppose that $q\in H_{\mathrm{loc}}^{-1}(\mathbb{R})$ and
$\lambda_{0}(q)>0$. Then the equation $L_{q}y=0$ has two linearly independent
positive solutions $y_{1},y_{2}\in H_{\mathrm{loc}}^{1}(\mathbb{R})$.
\end{proposition}

\begin{remark}
For potentials $q \in L_{\mathrm{loc}}^{1}(\mathbb{R})$, the result above is
due to Murata \cite{Murata:1986}, Remark after Theorem 2.7.
\end{remark}

In the proof of Proposition \ref{prop.2pos} we will use

\begin{lemma}
\label{L:piecewise} Assume that $y\in H^{1}_{\mathrm{loc}}({\mathbb{R}})$, and
there exists a discrete subset $S\subset{\mathbb{R}}$, such that $L_{q}y=0$ on
${\mathbb{R}}\setminus S$. Then
\begin{equation}
\label{E:piecewise}L_{q}y=\sum_{z\in S}(u(z-0)-u(z+0))\delta(\cdot-z),
\end{equation}
where $u=y^{\prime}-Qy$, $Q^{\prime}=q$ as in \eqref{eq.sys}. In other words,
\begin{equation}
\label{E:dpiecewise}(y, L_{q}\varphi)=\sum_{z\in S}(u(z-0)-u(z+0))\overline
{\varphi(z)},
\end{equation}
for every $\varphi\in\mathcal{C}_{0}^{\infty}(\mathbb{R})$.
\end{lemma}

\begin{proof}
Using partition of unity, we can split any function $\varphi\in\mathcal{C}%
_{0}^{\infty}({\mathbb{R}})$ into a finite sum of functions $\varphi_{k}$ such
that for every $k$ a neighborhood of $\mathrm{supp}\,\varphi_{k}$ contains at
most one point from $S$. Therefore, taking into account translation
invariance, we see that it suffices to consider the case when $S=\{0\}$. So we
will assume that $y\in H^{1}_{\mathrm{loc}}({\mathbb{R}})$ and $L_{q}y=0$ on
${\mathbb{R}}\setminus0$.

Integrating by parts (see \eqref{E:parts}) and using \eqref{eq.sys}, we get
\begin{align*}
(y,L_{q}\varphi) =(y^{\prime},\varphi^{\prime}) + (qy,\varphi) \hskip2in\\
=(u+Qy,\varphi^{\prime}) + (qy,\varphi) =(u,\varphi^{\prime}) +
(qy-(Qy)^{\prime},\varphi)\\
=(u,\varphi^{\prime}) - (Qy^{\prime},\varphi) =(u,\varphi^{\prime}) -
(Q^{2}y+Qu,\varphi).
\end{align*}
Integrating by parts in the first term in the right hand side we obtain
\[
(u,\varphi^{\prime})=\int_{-\infty}^{0} u\overline{\varphi^{\prime}}dx +
\int^{\infty}_{0} u\overline{\varphi^{\prime}}dx= (u(-0)-u(+0))\overline
{\varphi(0)} - ([u^{\prime}],\varphi),
\]
where $[u^{\prime}]$ is the locally integrable function on ${\mathbb{R}}$
which coincides with $u^{\prime}$ on ${\mathbb{R}}\setminus\{0\}$. Since
$[u^{\prime}]=-Q^{2}y-Qu$ due to \eqref{eq.sys}, we finally obtain
\[
(y,L_{q}\varphi)=(u(-0)-u(+0))\overline{\varphi(0)},
\]
as required.
\end{proof}

\begin{corollary}
\label{C:piecewise} Let $y$ satisfy the conditions of Lemma \ref{L:piecewise}
and have a compact support (so that $S$ can be taken finite). Then
\begin{equation}
\mathfrak{t}_{q}(y,y)=\int_{\mathbb{R}}(|y^{\prime}|^{2}+q|y|^{2}%
)dx=\sum_{z\in S}(u(z-0)-u(z+0))\overline{y(z)}. \label{E:qpiecewise}%
\end{equation}

\end{corollary}

\begin{proof}
Taking limit in \eqref{E:dpiecewise} over a sequence $\varphi_{k}$ converging
to $y$ in $H^{1}_{\mathrm{comp}}({\mathbb{R}})$, we obtain \eqref{E:qpiecewise}.
\end{proof}

\begin{proof}
[Proof of Proposition \ref{prop.2pos}]Using notations from the proof of
Proposition \ref{prop.pos} (see \eqref{E:yc}), for any $c>0$ define a test
function (to use in \eqref{eq.qlambda})
\begin{equation}
\label{E:psic}\psi_{c}(x)=
\begin{cases}
y_{-c}(x), & x\in(-c,0);\\
y_{c}(x), & x\in[0,c);\\
0, & x\not \in (-c,c).
\end{cases}
\end{equation}
Clearly, $\psi_{c}\in H^{1}_{\mathrm{comp}}(\mathbb{R})$, $\psi_{c}(0)=1$, and
$L_{q}\psi_{c}=0$ on $\mathbb{R}\setminus S$ where $S=\{-c,0,c\}$.

Applying Corollary \ref{corollary.dirichlet} to $y_{c}-y_{c^{\prime}}$, we see
that $c\mapsto\psi_{c}(x)$ is an increasing function of $c>0$ for any fixed
$x\in{\mathbb{R}}$. Therefore, the $L^{2}$-norm $\|\psi_{c}\|$ increases with
$c$ as well.

By Corollary \ref{C:piecewise} we have
\[
\mathfrak{t}_{q}(\psi_{c},\psi_{c})=u_{-c}(0)-u_{c}(0).
\]
It follows from Lemma \ref{lemma.zero} that $u_{-c}(0)$ decreases and
$u_{c}(0)$ increases as $c$ increases. Therefore, $c\mapsto\mathfrak{t}%
_{q}(\psi_{c},\psi_{c})$ is decreasing as $c$ increases. It follows that the
fraction in \eqref{eq.qlambda} is decreasing as well. To prove the desired
statement, it is enough to establish that the limit of this fraction is $0$ as
$c\to+\infty$, provided we know that the equation $L_{q}y=0$ has only one
positive solution with $y(0)=1$. To this end note that the limits of $y_{-c}$
and $y_{c}$ exists and are both positive solutions, according to the arguments
given in the proof of Proposition \ref{prop.pos}. Due to our uniqueness of
positive solution assumption these limits should coincide. But then we should
also have
\[
\lim_{c\to+\infty}u_{-c}(0)=\lim_{c\to+\infty}u_{c}(0),
\]
because $y_{-c}(0)=y_{c}(0)=1$ and the map $y\mapsto\{y(0),u(0)\}$ is a linear
topological isomorphism between the space of all solutions of $L_{q}y=0$ with
the $H^{1}_{\mathrm{loc}}(\mathbb{R})$-topology and the space $\mathbb{C}^{2}$
 It follows that
\[
\lim_{c\to+\infty}\mathfrak{t}_{q}(\psi_{c},\psi_{c})=0,
\]
which implies the desired statement.
\end{proof}

\begin{proof}
[Proof of Theorem \ref{thm.nofold}]To prove part (i), it is enough to show
that $E_{>}\subset E_{2}$ since it then follows by taking complements that
$E_{1}\subseteq E_{\bullet}$. In Proposition \ref{prop.2pos}, we established
that any $q\in\operatorname{Im}(B)$ with $\lambda_{0}(q)>0$ has two linearly
independent, positive solutions of $L_{q}y=0$ in $H_{\mathrm{loc}}%
^{1}(\mathbb{R})$, so $E_{>}\subset E_{2}$.

To prove part (ii), we first note that, by Lemma \ref{lemma.e2}, there are
compactly supported potentials in $E_{2}$, and by the remark above,
$\lambda_{0}(q)=0$ for such potentials, so $E_{2}\cap E_{\bullet}$ is
nonempty. Next, note that $B\left(  \mathcal{C}_{0}^{\infty}(\mathbb{R)}%
\right)  \subset\mathcal{C}_{0}^{\infty}(\mathbb{R})$ so $B\left(
\mathcal{C}_{0}^{\infty}(\mathbb{R)}\right)  \subset E_{\bullet}$. On the
other hand, by Lemma \ref{lemma.dense1}, $B\left(  \mathcal{C}_{0}^{\infty
}(\mathbb{R)}\right)  $ is dense in $\operatorname{Im}(B)$, so $E_{\bullet}$
is dense in $\operatorname{Im}(B)$. We have already shown that $E_{1}$ is
dense in $\operatorname{Im}(B)$, so $E_{1}$ is also dense in $E_{\bullet}$ by
part (i). The proof of Lemma \ref{lemma.dense2} shows that $E_{2}\cap
E_{\bullet}$ is dense in $E_{\bullet}.$
\end{proof}

\begin{remark}
Note that the map $\Phi:E_{\bullet}\times\mathbb{R}_{\geq0}\rightarrow
\operatorname{Im}(B)$ defined by $(q,c)\mapsto q+c$ is continuous and
bijective, but not a homeomorphism. Otherwise, $\Phi(E_{\bullet}\times\left\{
0\right\}  )\subset\operatorname{Im}(B)$ would be closed, and, as $E_{1}$ is
dense in $\operatorname{Im}(B)$, we conclude that $E_{\bullet}%
=\operatorname{Im}(B)$, a contradiction. The interpretation that $E_{2}$ is at
least \textquotedblleft one dimension larger\textquotedblright\ than $E_{1}$
could therefore be somewhat misleading. Note that the inverse of $\Phi$ is
given by $\Phi^{-1}:\operatorname{Im}(B)\rightarrow E_{\bullet}\times
\mathbb{R}_{\geq0}$, $q\mapsto\left(  q-\lambda_{0}(q),\lambda_{0}(q)\right)
$. Hence, $\Phi^{-1}$ not being continuous means that $q\mapsto\lambda_{0}(q)$
is not continuous in $H_{\mathrm{loc}}^{-1}(\mathbb{R})$.
\end{remark}

\section{Application to KdV}

\label{sec.KdV}

In this section we apply our results on the Miura map to prove existence of
solutions of the Korteweg-de Vries equation in $H^{-1}(\mathbb{R})$ for
initial data in the range $\operatorname{Im}(B_{0})$ of the Miura map
$B_{0}:L^{2}(\mathbb{R})\rightarrow H^{-1}(\mathbb{R})$. We follow the
approach of Tsutsumi \cite{Tsu:1989}, \ who proved such an existence result
for initial data a positive, finite Radon measure on $\mathbb{R}$. His
arguments combined with our results on the Miura map $B_{0}$ lead to the
following theorem. Recall that, for a real-valued distribution $u\in
H^{-1}(\mathbb{R})$, $\left[  u\right]  $ denotes the special integral of $u$
(see Theorem \ref{thm.special} and the discussion that precedes it).

\begin{theorem}
\label{thm.6.1}Assume that $u_{0}\in\operatorname{Im}(B_{0})$. Then there
exists a global weak solution of KdV with $u(t)\in\operatorname{Im}\left(
B_{0}\right)  $ for all $t\in\mathbb{R}$. More precisely: \newline(i) $u\in
L^{\infty}(\mathbb{R},H^{-1}(\mathbb{R}))\cap L_{\mathrm{loc}}^{2}%
(\mathbb{R}^{2})$, \newline(ii) for all functions $\varphi\in\mathcal{C}%
_{0}^{\infty}(\mathbb{R}^{2})$, the identity%
\[
\int_{\mathbb{R}}\int_{\mathbb{R}}\left(  -u\varphi_{t}-u\varphi_{xxx}%
+3u^{2}\varphi_{x}\right)  ~dx~dt=0
\]
holds, \newline(iii) $\lim_{t\rightarrow0}u(t)=u_{0}$ in $H^{-1}(\mathbb{R})$,
and\newline(iv) $0\leq\left[  u(t)\right]  \leq\left[  u_{0}\right]  $ for all
$t\in\mathbb{R}$ and $\lim_{t\rightarrow0} \left[  u(t) \right]  = \left[
u_{0}\right]  .$
\end{theorem}

\begin{remark}
\ Recall that $u_{0} \in\operatorname{Im}(B_{0})$ means that $u_{0}\in
H^{-1}(\mathbb{R})$ with the property that $L_{u_{0}} \geq0 $ and $\left[
u_{0}\right]  $ exists. Instead of the assumption for $\left[  u_{0}\right]  $
to exist, one can equivalently assume that \thinspace$u_{0}=f^{\prime}+g$ for
some functions $f\in L^{2}(\mathbb{R})$ and $g\in L^{1}(\mathbb{R})$ -- see
Theorem \ref{thm.im2} and Proposition \ref{prop.image2}.
\end{remark}

To prove Theorem \ref{thm.6.1}, we need to recall a result of Kato
\cite{Kato:1983} and, independently, of Kruzhkov and Faminski\u{\i}
\cite{KF:1984} -- see also \cite{BS:1993} and \cite{GTV:1990}. Consider the
modified Korteweg-de Vries equation (mKdV)%
\begin{equation}
\partial_{t}v=-\partial_{x}^{3}v+6v^{2}\partial_{x}v \label{eq.6.1}%
\end{equation}
with initial data%
\begin{equation}
v(0)=v_{0}. \label{eq.6.2}%
\end{equation}

\begin{theorem}
\label{thm.6.2} \cite{Kato:1983}, \cite{KF:1984} Let $v_{0}\in L^{2}%
(\mathbb{R})$. Then the initial value problem (\ref{eq.6.1})-(\ref{eq.6.2})
has a weak solution in $L^{\infty}(\mathbb{R},L^{2}(\mathbb{R}))$. More
precisely, $v$ satisfies:\newline(i) $v\in L^{\infty}(\mathbb{R}%
,L^{2}(\mathbb{R}))\cap L_{\mathrm{loc}}^{2}(\mathbb{R},H_{\mathrm{loc}}%
^{1}(\mathbb{R})\mathbb{)}$,\newline(ii) the identity%
\[
\int_{\mathbb{R}}~\int_{\mathbb{R}}\left(  -v\varphi_{t}-v\varphi_{xxx}%
+2v^{3}\varphi_{x}\right)  ~dt~dx=0
\]
holds for all $\varphi\in\mathcal{C}_{0}^{\infty}(\mathbb{R}^{2})$%
,\newline(iii) $\lim_{t\rightarrow0}v(t)=v_{0}$ in $L^{2}(\mathbb{R})$,
and\newline(iv) $\left\Vert v(t)\right\Vert _{L^{2}(\mathbb{R})}\leq\left\Vert
v_{0}\right\Vert _{L^{2}(\mathbb{R})}$ for all $t\in\mathbb{R}$.
\end{theorem}

The following result improves the one of Tsutsumi \cite{Tsu:1989} by adapting it to our more
general setting, and relies on the identity (\ref{eq.MiuraId}).

\begin{proposition}
\label{prop.6.3} Let $v=v(t)$ be a solution of (\ref{eq.6.1})-(\ref{eq.6.2})
with $v_{0}\in L^{2}(\mathbb{R})$ and the properties listed in Theorem
\ref{thm.6.2}. Let $u_{0}:=v_{0}^{\prime}+v^{2}_{0}$. Then $u:=v^{\prime
}+v^{2}$ is a solution of KdV. More precisely:\newline(i) $u\in L^{\infty
}(\mathbb{R},H^{-1}(\mathbb{R}))\cap L_{\mathrm{loc}}^{2}(\mathbb{R}^{2})$,
\newline(ii) the identity%
\[
\int_{\mathbb{R}}~\int_{\mathbb{R}}\left(  -u\varphi_{t}-u\varphi_{xxx}%
+3u^{2}\varphi_{x}\right)  ~dt~dx=0
\]
holds for all $\varphi\in\mathcal{C}_{0}^{\infty}(\mathbb{R}^{2})$,
\newline(iii) $\lim_{t\rightarrow0}u(t)=u_{0}$ in $H^{-1}(\mathbb{R})$,
\newline(iv) $0\leq\left[  u(t)\right]  \leq\left[  u_{0}\right]  $ for all
$t\in\mathbb{R},$ and $\lim_{t\rightarrow0} \left[  u(t) \right]  = \left[
u_{0}\right]  .$
\end{proposition}

\begin{proof}
Statement (i)\ follows from Theorem \ref{thm.6.2}(i) together with the fact
that $B_{0}:L^{2}(\mathbb{R})\rightarrow H^{-1}(\mathbb{R})$ is a bounded,
continuous map -- see Proposition \ref{prop.Miura.continuous}. Statement (iii)
follows from Theorem \ref{thm.6.2}(iii) and the continuity of $B_{0}$ whereas
the claimed inequality in (iv)\ follows from Theorem \ref{thm.6.2}(iv) and the
fact that $\left[  u(t)\right]  =\left\Vert v(t)\right\Vert ^{2}_{L^{2}(\mathbb{R})}$ -- see
formula (\ref{eq.r}). To prove the second statement in (iv), note that by
Theorem \ref{thm.6.2}(iii), $\lim_{t\rightarrow0} \left\Vert v(t)\right\Vert
_{L^{2}(\mathbb{R})} = \left\Vert v_{0} \right\Vert _{L^{2}(\mathbb{R})}$. As
$\left[  u(t)\right]  =\left\Vert v(t)\right\Vert ^{2} _{L^{2}(\mathbb{R})},$
the second statement in (iv) then follows as well. Statement (ii)\ is proved
in \cite{Tsu:1989}. For the convenience of the reader we include a detailed
proof of it.

Let $\rho:\mathbb{R}\times\mathbb{R}\rightarrow\mathbb{R}$ be the smooth
mollifier, i.e. a smooth positive function with support in the unit disc in
$\mathbb{R }^{2}$, $\rho(0,0) > 0,$ and normalized by $\int_{\mathbb{R}^{2}%
}\rho(t,x)~dt~dx=1$. For $\varepsilon>0$, set
\[
\rho_{\varepsilon}\left(  t,x\right)  :=\frac{1}{\varepsilon^{2}}\rho\left(
\frac{t}{\varepsilon},\frac{x}{\varepsilon}\right)  .
\]

Given the solution $v(t)$ of mKdV, define for $\left(  t,x\right)
\in\mathbb{R}^{2}$%
\begin{align*}
v_{\varepsilon}(t,x)  &  :=\left(  \rho_{\varepsilon}\ast v\right)  (t,x)\\
&  =\int_{\mathbb{R}^{2}}\rho\left(  t-s,x-y\right)  ~v\left(  s,y\right)
~ds~dy.
\end{align*}
Note that for any $\varepsilon>0$ and $\left(  t,x\right)  \in\mathbb{R}^{2}$,
$\rho_{\varepsilon}\left(  t-s,x-y\right)  \in\mathcal{C}_{0}^{\infty
}(\mathbb{R}^{2})$ as a function of $\left(  s,y\right)  \in\mathbb{R}^{2}$.
Further define the function $u_{\varepsilon} \in\mathcal{C}_{0}^{\infty
}(\mathbb{R}^{2}),$
\[
u_{\varepsilon}:=\frac{\partial}{\partial x}v_{\varepsilon}+v_{\varepsilon
}^{2}.
\]
According to (\ref{eq.MiuraId}), one has%
\begin{align}
\frac{\partial}{\partial t}u_{\varepsilon}+\frac{\partial^{3}}{\partial x^{3}%
}u_{\varepsilon}-6u_{\varepsilon}\frac{\partial}{\partial x}u_{\varepsilon}
&  =\left(  \frac{\partial}{\partial x}+2v_{\varepsilon}\right)  \left(
\frac{\partial}{\partial t}v_{\varepsilon}+\frac{\partial^{3}}{\partial x^{3}%
}v_{\varepsilon}-6\rho_{\varepsilon}\ast\left(  v^{2}\frac{\partial}{\partial
x}v\right)  \right) \label{eq.6.3}\\
&  +6\left(  \frac{\partial}{\partial x}+2v_{\varepsilon}\right)  \left(
\rho_{\varepsilon}\ast\left(  v^{2}\frac{\partial}{\partial x}v\right)
-v_{\varepsilon}^{2}\frac{\partial}{\partial x}v_{\varepsilon}\right)
.\nonumber
\end{align}
By assumption, \ $v$ is a weak solution of mKdV, hence%
\begin{align*}
\frac{\partial}{\partial t}v_{\varepsilon}+\frac{\partial^{3}}{\partial x^{3}%
}v_{\varepsilon}-6\rho_{\varepsilon}\ast\left(  v^{2}\frac{\partial}{\partial
x}v\right)   &  =\rho_{\varepsilon}\ast\left(  \frac{\partial}{\partial
t}v+\frac{\partial^{3}}{\partial x^{3}}v-6v^{2}\frac{\partial}{\partial
x}v\right) \\
&  =0.
\end{align*}
Multiplying (\ref{eq.6.3}) by an arbitrary test function $\varphi
\in\mathcal{C}_{0}^{\infty}(\mathbb{R}^{2})$ and integrating by parts, one
obtains%
\begin{align}
&  \int_{\mathbb{R}^{2}}\left(  -u_{\varepsilon}\frac{\partial}{\partial
t}\varphi-u_{\varepsilon}\frac{\partial^{3}}{\partial x^{3}}\varphi
+3u_{\varepsilon}^{2}\frac{\partial}{\partial x}\varphi\right)
~dt~dx\label{eq.6.4}\\
&  =6\int_{\mathbb{R}^{2}}\left[  \rho_{\varepsilon}\ast\left(  v^{2}
\frac{\partial}{\partial x}v\right)  -v_{\varepsilon}^{2}\frac{\partial
}{\partial x}v_{\varepsilon}\right]  \left(  -\frac{\partial}{\partial
x}\varphi+2v_{\varepsilon}\varphi\right)  ~dt~dx\nonumber\\
~~~  &  =2\int_{\mathbb{R}^{2}}\left(  \rho_{\varepsilon}\ast v^{3}%
-v_{\varepsilon}^{3}\right)  \left(  \frac{\partial^{2}}{\partial x^{2}%
}\varphi-2\varphi\frac{\partial}{\partial x}v_{\varepsilon}-2v_{\varepsilon
}\frac{\partial}{\partial x}\varphi\right)  ~dt~dx.\nonumber
\end{align}
By Theorem \ref{thm.6.2}(i),
\[
\frac{\partial}{\partial x}v_{\varepsilon}\rightarrow\frac{\partial}{\partial
x}v~\text{in }L_{\mathrm{loc}}^{2}(\mathbb{R}^{2}).
\]
Lemma \ref{lemma.6bis} below together with Theorem \ref{thm.6.2}(i) implies
that $v\in L_{\mathrm{loc}}^{6}(\mathbb{R}^{2})$. Hence, \
\[
v_{\varepsilon}^{2}\rightarrow v^{2},~~v_{\varepsilon}^{3}\rightarrow
v^{3}~~\text{in }L_{\mathrm{loc}}^{2}(\mathbb{R}^{2}).
\]
Combining all of this, one obtains%
\begin{align*}
\rho_{\varepsilon}\ast v^{3}-v^{3}_{\varepsilon}  &  =\left(  \rho
_{\varepsilon}\ast v^{3}-v^{3}\right)  - \left(  v_{\varepsilon}^{3}%
-v^{3}\right)
  \rightarrow0 \,\text{  as  }\,\varepsilon\rightarrow0 \,\text{  in  }\, L_{\mathrm{loc}%
}^{2}(\mathbb{R}^{2})
\end{align*}
and%
\[
2\varphi\frac{\partial}{\partial x}v_{\varepsilon} + 2v_{\varepsilon}
\frac{\partial}{\partial x}\varphi\rightarrow2\varphi\frac{\partial}{\partial
x}v + 2v\frac{\partial}{\partial x}\varphi\, \text{   in  }\, L^{2}(\mathbb{R}^{2}).
\]
Since $u_{\varepsilon}=\partial_{x}v_{\varepsilon}+v_{\varepsilon}^{2}$, \ it
follows that%
\[
u_{\varepsilon}\rightarrow u\text{ in }L_{\mathrm{loc}}^{2}(\mathbb{R}^{2})
\]
and, as $\operatorname*{supp}\varphi$ is compact,
\[
\int_{\mathbb{R}^{2}}\left(  \rho_{\varepsilon}\ast v^{3}-v_{\varepsilon}%
^{3}\right)  \left(  \frac{\partial^{2}}{\partial x^{2}}\varphi - 2\varphi
\frac{\partial}{\partial x}v_{\varepsilon} - 2v_{\varepsilon}\frac{\partial
}{\partial x}\varphi\right)  ~dt~dx\rightarrow0
\]
as $\varepsilon\rightarrow0$. Therefore, taking the limit $\varepsilon
\rightarrow0$ in (\ref{eq.6.4}), we conclude that%
\[
\int_{\mathbb{R}^{2}}\left(  -u\frac{\partial}{\partial t}\varphi
-u\frac{\partial^{3}}{\partial x^{3}}\varphi+3u^{2}\frac{\partial}{\partial
x}\varphi\right)  ~dt\;dx=0.
\]

\end{proof}

\begin{lemma}
\label{lemma.6bis}Let $Q:=I\times J$ with $I:=\left[  -T,T\right]  $ and
$J:=\left[  -R,R\right]  $ where $T>0$ and $R>0$. Let
\[
\mathcal{B}_{I,J}=L^{2}\left(  I,H^{1}(J)\right)  \cap L^{\infty}(I,L^{2}(J))
\]
with norm%
\[
\left.  \operatorname*{ess}\sup\right.  _{t\in I}\left\Vert f(t,~\cdot
~)\right\Vert _{L^{2}(J)}+\left(  \int_{I}\left\Vert f(t,~\cdot~)\right\Vert
_{H^{1}(J)}^{2}~dt\right)  ^{1/2}.
\]
Then, \ for any $f\in\mathcal{B}_{I,J}$, \ the inequality%
\begin{equation}
\left\Vert f\right\Vert _{L^{6}\left(  Q\right)  }^{6}\leq C\left.
\operatorname*{ess}\sup\right.  _{t\in I}\left\Vert f(t,~\cdot~)\right\Vert
_{L^{2}(J)}^{4} \cdot\int_{I}\left\Vert f(t,~\cdot~)\right\Vert _{H^{1}%
(J)}^{2}~dt \label{eq.6.5}%
\end{equation}
holds, where $C>0$ is a constant which depends only on $R$. In particular,
$\mathcal{B}_{I,J}$ embeds continuously into $L^{6}(Q)$.
\end{lemma}

\begin{proof}
Let us assume first that $f\in\mathcal{C}_{0}^{\infty}(\mathbb{R}^{2})$.  
By the standard Sobolev embedding theorem, the space $H^{1/3}(J)$ is densely and continuously
embedded in $L^{6}(J)$ (see e.g. \cite{RS:1996}, Theorem 1, p 82),
\[
\left\Vert g\right\Vert _{L^{6}(J)}\leq C_{J}^{\prime}\left\Vert g\right\Vert
_{H^{1/3}(J)}%
\]
for all $g\in H^{1/3}(J)$ and some positive constant $C_{J}^{\prime}$.
Further, by interpolation, one has for any $g\in H^{1}(J)$ (see e.g.
\cite{RS:1996}, Remark 2, p 87)%
\[
\left\Vert g\right\Vert _{H^{1/3}(J)}\leq C_{J}^{\prime\prime}\left\Vert
g\right\Vert _{H^{1}(J)}^{1/3}\left\Vert g\right\Vert _{L^{2}(J)}^{2/3}.
\]
Setting $C=\left(  C_{J}^{\prime}C_{J}^{\prime\prime}\right)  ^{6}$, we
see that 
\begin{align*}
\int_{I\times J}\left\vert f(t,x)\right\vert ^{6} ~dt~dx  &  \leq\int
_{I}\left(  C_{J}^{\prime}\left\Vert f(t,~\cdot~)\right\Vert _{H^{1/3}%
(J)}\right)  ^{6}~dt\\
&  \leq C\int_{I}\left\Vert f\left(  t,~\cdot~\right)  \right\Vert _{H^{1}%
(J)}^{2}\left\Vert f\left(  t,~\cdot~\right)  \right\Vert _{L^{2}(J)}^{4}~dt.
\end{align*}
By approximation, the above inequality holds for any $f \in \mathcal{B}_{I,J}.$ For
such $f$ we have
\begin{align*}
&\int_{I}\left\Vert f\left(  t,~\cdot~\right)  \right\Vert _{H^{1}%
(J)}^{2}\left\Vert f\left(  t,~\cdot~\right)  \right\Vert _{L^{2}(J)}^{4}~dt\\
&  \leq \left.  \operatorname*{ess}\sup\right.  _{t\in I}\left\Vert f\left(
t,~\cdot~\right)  \right\Vert _{L^{2}(J)}^{4}\int_{I}\left\Vert f\left(
t,~\cdot~\right)  \right\Vert _{H^{1}(J)}^{2}~dt.
\end{align*}
Combinig the two previous inequalities ends the proof.
\end{proof}

\begin{proof}
[Proof of Theorem \ref{thm.6.1}]The proof is the one given in Tsutsumi
\cite{Tsu:1989}, adapted to our more general setting. By our assumption, 
$B_{0}(v_{0})=u_{0}$ for some $v_{0}\in L^{2}(\mathbb{R})$. 
By Theorem \ref{thm.6.2}, there exists a solution $v\in
L^{\infty}\left(  \mathbb{R},L^{2}(\mathbb{R})\right)  \cap L_{\mathrm{loc}%
}^{2}\left(  \mathbb{R},H_{\mathrm{loc}}^{1}(\mathbb{R})\right)  $ of
(\ref{eq.6.1})-(\ref{eq.6.2}). By Proposition \ref{prop.6.3}, $u(t):=B(v(t))$
is a solution of KdV\ satisfying (i)-(iv).
\end{proof}

\appendix{}

\section{Positive Solutions for Square-Well Potentials}

\label{sec.square}

In this appendix we present some elementary but important examples of
potentials in $\operatorname{Im}(B)$ and the associated positive solutions.
Let
\[
q_{a,b}(x)=\left\{
\begin{array}
[c]{cc}%
b^{2}, & -a<x<a\\
& \\
0, & \left\vert x\right\vert \geq a
\end{array}
\right.
\]
where $a,b >0.$ It is easy to see that%
\begin{equation}
y_{+}(x)=\left\{
\begin{array}
[c]{lll}%
1/\cosh(ba), &  & x<-a\\
&  & \\
\cosh(b(x+a))/\cosh(ba), &  & -a<x<a\\
&  & \\
\left(\cosh(2ab)+b(x-a)\sinh(2ab)\right)/\cosh(ba), &  & x>a
\end{array}
\right.  \label{eq.square}%
\end{equation}
and $y_{-}(x):=y_{+}(-x)$ are linearly independent positive solutions
of $-y^{\prime\prime}+q_{a,b}y=0$.

If $\lambda>0$ and $b=\left(  \lambda/2a\right)  ^{1/2}$ then $\int
q_{a,b}(x)~dx=\lambda$. Taking $a\downarrow0$ we recover in the limit
$q=\lambda\delta$ where $\delta$ is the Dirac $\delta$-distribution at $x=0$.
In this limit%
\[
y_{+}(x)=\left\{
\begin{array}
[c]{lll}%
1 &  & x<0\\
&  & \\
1+\lambda x &  & x>0
\end{array}
\right.
\]
and again $y_{-}(x)=y_{+}(-x)$.

Let us determine the preimage of $q = \lambda\delta$ by the Miura map. From
the explicit formulas we have $\int_{0}^{\infty}y_{+}(s)^{-2}ds<\infty$ but
$\int_{-\infty}^{0}y_{+}(s)^{-2}ds=+\infty$, while the reverse is true for
$y_{-}$. If $H$ is the Heaviside function
\[
H(x)=\left\{
\begin{array}
[c]{cc}%
0 & x<0\\
& \\
1 & x>0
\end{array}
\right.
\]
then the logarithmic derivatives%
\[
\frac{y_{+}^{\prime}(x)}{y_{+}(x)}=\frac{\lambda H(x)}{1+\lambda x}%
,~~\frac{y_{-}^{\prime}(x)}{y_{-}(x)}=-\frac{y_{+}^{\prime}(-x)}{y_{+}(-x)}%
\]
belong to $L^{2}(\mathbb{R})$. Hence
\[
B^{-1}(\lambda\delta) = \left\{ (1 - \theta) \frac{\lambda H(x)}{1+ \lambda x} -
\theta\frac{\lambda H(-x)} {1 - \lambda x} \mid 0 \le\theta\le1 \ \right\}.
\]
\newline

\section{Positive Schr\"odinger Operators}

\label{sec.green}

\noindent In this Appendix we provide more information about 
 Schr\"{o}dinger operators $L_{q}$ which are positive or, more
generally, semibounded below, and
have real potentials $q\in H_{\mathrm{loc}}^{-1}(\mathbb{R})$. Namely, we will
show that the corresponding quadratic form (defined on $\mathcal{C}%
_{0}^{\infty}(\mathbb{R})$) is closable and describe the domain of its
closure. We will also describe the domain of the corresponding self-adjoint
operator. Finally, for strictly positive $L_{q}$ (such that $\lambda_{0}%
(q)>0$, see \eqref{eq.lambda}, \eqref{eq.qlambda}) we construct Green's
function and use it to give an alternative proof of Proposition
\ref{prop.2pos}.

The case of semi-bounded $L_{q}$ for many purposes is reduced to the case when
$L_{q}\geq0$ or even to the case when $L_{q}$ is strictly positive (that is
$\lambda_{0}(q)>0$) by adding a sufficiently large constant to $q$. So let us
assume first that $L_{q}\geq0$. By Theorem \ref{thm.im1} there exists a
function $r\in L_{\mathrm{loc}}^{2}(\mathbb{R})$ such that $q=B(r)$. Then
$L_{q}$ admits a formal factorization \eqref{E:factor}, i.e., a presentation
$L_{q}=P^{+}P$, where $P=\left(  \partial_{x}-r\right)  $ and $P^{+}=-\left(
\partial_{x}+r\right)  $, so that $P^{+}$ is the operator formally adjoint to
$P$ in $L^{2}(\mathbb{R})$.

Clearly, $P,P^{+}$ are well defined on the space $\mathcal{C}_{0}^{\infty
}(\mathbb{R})$ which is dense in $L^{2}(\mathbb{R})$, so that $P,P^{+}$ map
$\mathcal{C}_{0}^{\infty}(\mathbb{R})$ to $L^{2}(\mathbb{R})$ and
\[
(Pu,v)=(u,P^{+} v), \quad u,v\in\mathcal{C}_{0}^{\infty}(\mathbb{R}).
\]
It follows that the operators $P,P^{+}$ are closable, with the closures which
we will denote by $\overline{P}, \overline{P^{+}}$. They also have
adjoint operators in $L^{2}(\mathbb{R})$, which will be denoted $P^{*},
(P^{+})^{*}$. These operators are closed extensions of $P^{+}, P$
respectively. Since $P^{*}, (P^{+})^{*}$ are closed, we have
\begin{equation}
\label{E:adjinclusions}\overline{P}\subset(P^{+})^{*}, \quad\overline{P^{+}%
}\subset P^{*}.
\end{equation}

\begin{lemma}
\label{operatorP} (i) We have
\begin{equation}
\label{E:adjoints}\overline{P}=(P^{+})^{*}, \quad\overline{P^{+}}=P^{*}.
\end{equation}

(ii) The domains of the operators in \eqref{E:adjoints} are as follows:
\begin{equation}
\label{E:DP}\mathfrak{D}(\overline{P})= \left\{  u\in L^{2}(\mathbb{R})\cap
W_{\mathrm{loc}}^{1,1}(\mathbb{R}): Pu\in L^{2}(\mathbb{R})\right\}  ,
\end{equation}
\begin{equation}
\label{E:DPstar}\mathfrak{D}(\overline{P^{+}})= \left\{  v\in L^{2}%
(\mathbb{R})\cap W_{\mathrm{loc}}^{1,1}(\mathbb{R}): P^{+}v\in L^{2}%
(\mathbb{R})\right\}  ,
\end{equation}
where the operators $P,P^{+}$ are applied in the usual distributional sense.

(iii) $\mathcal{C}_{0}^{\infty}(\mathbb{R})$ is an operator core for each of
the operators in \eqref{E:adjoints}.
\end{lemma}

\begin{proof}
It is easy to see that the right-hand sides in \eqref{E:DP}, \eqref{E:DPstar}
coincide with the domains of the adjoint operators $(P^{+})^{\ast},P^{\ast}$
respectively. Indeed, the relation $(P^{+})^{\ast}u=f$ means that $u,f\in
L^{2}(\mathbb{R})$ and for every $\varphi\in\mathcal{C}_{0}^{\infty
}(\mathbb{R})$
\[
(u,P^{+}\varphi)=-(u,\partial_{x}\varphi)-(u,r\varphi)=(f,\varphi).
\]
Since $ru\in L_{\mathrm{loc}}^{1}(\mathbb{R})$, this is equivalent to
$\partial_{x}u-ru=f$, where $\partial_{x}$ is applied in the sense of
distributions. It follows that $\partial_{x}u=f+ru\in L_{\mathrm{loc}}%
^{1}(\mathbb{R})$, hence $u\in W_{\mathrm{loc}}^{1,1}(\mathbb{R})$ or,
equivalently, $u$ is absolutely continuous. This means that the right hand
side of \eqref{E:DP} coincides with $\mathfrak{D}((P^{+})^{\ast})$. The same
argument applies to the operator $P^{\ast}$ and the right-hand side of \eqref{E:DPstar}.

Taking into account the inclusions \eqref{E:adjinclusions}, we see that to
establish all statements of the lemma, it suffices to show that the right hand
sides of \eqref{E:DP}, \eqref{E:DPstar} belong to the domains of $\overline
{P}, \overline{P^{+}}$ respectively. This is easily done by use of 
Friedrichs' mollifiers. It is essentially a special case of Friedrichs'
\cite{Friedrichs:1944} well-known result on equality of weak and strong
extensions of differential operators, but we give the proof for the reader's
convenience. We will give the arguments for $P$ (the arguments for $P^{+}$ are
the same).

So let us assume that
\begin{equation}
u\in L^{2}(\mathbb{R})\cap W_{\mathrm{loc}}^{1,1}(\mathbb{R}),Pu\in
L^{2}(\mathbb{R}). \label{E:uin}%
\end{equation}
We need to show that $u$ may be approximated by a sequence $\left\{
u_{n}\right\}  _{n\geq1}$ from $\mathcal{C}_{0}^{\infty}(\mathbb{R})$ with
$u_{n}\rightarrow u$ and $Pu_{n}\rightarrow Pu$ in $L^{2}(\mathbb{R})$. First,
we show that it suffices to consider $u$ satisfying \eqref{E:uin} and
additionally having compact support. To this end, take $\chi\in\mathcal{C}%
_{0}^{\infty}(\mathbb{R})$. Then $\chi u$ also satisfies \eqref{E:uin} and
$P\left(  \chi u\right)  =\chi^{\prime}u+\chi Pu$. If $\chi_{n}\in
\mathcal{C}_{0}^{\infty}(\mathbb{R})$ satisfies the conditions $0\leq\chi
_{n}\leq1$, $\chi_{n}(x)=1$ for $\left\vert x\right\vert \leq n$, $\chi
_{n}(x)=0$ for $\left\vert x\right\vert \geq n+1$, and $\left\vert \chi
_{n}^{\prime}(x)\right\vert \leq2$, then $\chi_{n}u\rightarrow u$ in
$L^{2}(\mathbb{R})$ as $n\rightarrow\infty$. Moreover, $\chi_{n}^{\prime
}u\rightarrow0$ and $\chi_{n}Pu\rightarrow Pu$ in $L^{2}(\mathbb{R})$ as
$n\rightarrow\infty$, so $P(\chi_{n}u)=\chi_{n}^{\prime}u+\chi_{n}%
Pu\rightarrow Pu$ in $L^{2}(\mathbb{R})$ as $n\rightarrow\infty$ Thus, we may
assume that $u$ has compact support.

Given $u$ satisfying \eqref{E:uin} and having compact support, we now use
Friedrichs mollifiers to construct a sequence of approximants from
$\mathcal{C}_{0}^{\infty}(\mathbb{R})$. Let $j\in\mathcal{C}_{0}^{\infty
}(\mathbb{R})$ be a nonnegative function with $\int j(x)~dx=1$, and, for any
$k\in\mathbb{N}$, let $j_{k}(x)=kj(kx)$, and let $u_{k}=u\ast j_{k}$. Clearly
$u_{k}\in\mathcal{C}_{0}^{\infty}(\mathbb{R})$ and $u_{k}\rightarrow u$ in
$L^{2}(\mathbb{R})$. We claim that $Pu_{k}\rightarrow Pu$ in $L^{2}%
(\mathbb{R})$. Since $u\in W_{\mathrm{comp}}^{1,1}(\mathbb{R})$ and $r\in
L_{\mathrm{loc}}^{2}(\mathbb{R})$, it follows that $ru\in L^{1}(\mathbb{R})$.
Moreover, as $u$ satisfies \eqref{E:uin} and the support of $u$ is compact,
$Pu\in L^{1}(\mathbb{R}).$ Therefore $u^{\prime}=Pu+ru$ belongs to
$L^{1}(\mathbb{R})$. Hence, $u$ is a bounded, continuous function which shows
that $ru\in L^{2}(\mathbb{R}).$ As a consequence, $u^{\prime}=Pu+ru\in
L^{2}(\mathbb{R})$. Thus $u_{k}^{\prime}\rightarrow u^{\prime}$ in~$L^{2}%
(\mathbb{R}),$ $u_{k}\rightarrow u$ in $L^{\infty}(\mathbb{R})$ as
$k\rightarrow\infty$ and so%
\[
Pu_{k}-Pu=\left(  u_{k}^{\prime}-u^{\prime}\right)  +r\left(  u_{k}-u\right)
\]
converges to zero in $L^{2}(\mathbb{R})$ as $k\rightarrow\infty$.
\end{proof}

Now let us recall a classical theorem of von Neumann \cite{vonNeumann:1932}
(see also \cite{RS:1975}, Theorem XI.23) which asserts that if $A$ is a closed
densely defined operator in a Hilbert space, then the (generally unbounded)
operator $H=A^{\ast}A$ is self-adjoint. Here the domain of $A^{\ast}A$ is
naturally defined as
\[
\mathfrak{D}(A^{\ast}A)=\{u\in\mathfrak{D}(A),Au\in\mathfrak{D}(A^{\ast})\}.
\]
(Note that an essentially inverse statement also holds: if two densely defined
operators $A,A^{+}$ are formally adjoint, that is
\[
(Au,v)=(u,A^{+}v),\quad u\in\mathfrak{D}(A),v\in\mathfrak{D}(A^{+}),
\]
and $A^{+}A$ is essentially self-adjoint, then the closures $\overline
{A},\overline{A^{+}}$ are adjoint to each other; see the appendix to
\cite{Shubin:1976}.)

\medskip The following lemma is well-known.

\begin{lemma}
\label{L:AstarA} Let $A$ be a closed densely defined operator in a Hilbert
space, and $H=A^{*}A$. Denote by $\mathfrak{t}_{H}$ the quadratic form of $H$,
and let $\mathfrak{D}(\mathfrak{t}_{H})$ be its domain i.e. $\mathfrak{D}%
(\mathfrak{t}_{H})=\mathfrak{D}(H^{1/2})$. Then $\mathfrak{D}(\mathfrak{t}%
_{H})=\mathfrak{D}(A)$ and
\[
\mathfrak{t}_{H}(u,u)=\|Au\|^{2},\quad u\in\mathfrak{D}(A).
\]

\end{lemma}

\begin{proof}
Take the polar decomposition $A=U|A|$, where $|A|=(A^{*}A)^{1/2}=H^{1/2}$, and
$U$ partial isometry with $\mathrm{Ker}\, U=\mathrm{Ker}\, A$ (see e.g. Sect.
VIII.9 in \cite{RS:1975}). It remains to notice that $\mathfrak{D}%
(A)=\mathfrak{D}(|A|)$ (because $U$ is bounded), and
\[
\mathfrak{t}_{H}(u,u)=\|H^{1/2}u\|^{2}=\||A|u\|^{2}=\|Au\|^{2}, \quad
u\in\mathfrak{D}(A),
\]
because $U$ is an isometry on the range of $|A|$.
\end{proof}

Note that a positive self-adjoint operator $H$ is uniquely defined by its
(positive, closed) quadratic form (see e.g. Theorem VIII.15 in \cite{RS:1975}).

Taking the quadratic form $\mathfrak{t}_{q}$, corresponding to a potential
$q\in H_{\mathrm{loc}}^{-1}(\mathbb{R})$ and assuming that it is positive (or,
more generally, semi-bounded below), we can construct a unique self-adjoint
operator $H$ with this form. It follows from the considerations above that
when the form $\mathfrak{t}_{q}$ is positive, we can write this operator in
the form $H=P^{\ast}\overline{P}$, and the domain of $H$ is
\begin{equation}
\mathfrak{D}(H)=\{u\in L^{2}(\mathbb{R}),\;(\partial_{x}-r)u\in L^{2}%
(\mathbb{R}),\ (\partial_{x}+r)\left[  (\partial_{x}-r)u\right]  \in
L^{2}(\mathbb{R})\}, \label{E:domainLq}%
\end{equation}
where $r\in L_{\mathrm{loc}}^{2}(\mathbb{R})$ is a solution of the Riccati
equation $r^{\prime}+r^{2}=q$.

In case when $L_{q}$ is semibounded below but not positive, the arguments
given above should be applied to the operator $L_{q}+c$ with $c>0$ such that
$L_{q}+c\ge0$, with subsequent subtracting of the same constant $c$ from the
resulting operator (which does not change the domain of the operator). The
resulting operator $H$ will not depend of the choice of $c$ because different
choices of $c$ lead to the same (closed) quadratic form of the resulting operator.

The following lemma simplifies calculation of $Hu$ if we know that
$u\in\mathfrak{D}(H)$.

\begin{lemma}
\label{L:domain-exp} Let $q\in H^{-1}_{\mathrm{loc}}(\mathbb{R})$ be such that
$L_{q}\ge0$. Then $H$ can be extended to a linear operator $\widetilde{L}_{q}$
with the domain
\begin{equation}
\label{E:domain-exp}\mathfrak{D}(\widetilde{L}_{q})= \{u\in H^{1}%
_{\mathrm{loc}}(\mathbb{R})\cap L^{2}(\mathbb{R}),\; \widetilde{L}_{q}u\in
L^{2}(\mathbb{R})\},
\end{equation}
where $\widetilde{L}_{q}$ is $-\partial_{x}^{2}+q$ applied in the sense of
distributions, i.e. both $\partial_{x}^{2}$ and $q$ act as linear continuous
operators $H^{1}_{\mathrm{loc}}(\mathbb{R})\to H^{-1}_{\mathrm{loc}%
}(\mathbb{R})$ ($\partial_{x}^{2}$ acts as the distributional derivative, and
$q$ acts as a multiplier in these spaces).
\end{lemma}

\begin{proof}
If $u\in\mathfrak{D}(H)$ (as described by \eqref{E:domainLq}) then $u\in
L^{2}(\mathbb{R})$ and $f:=u^{\prime}-ru\in L^{2}(\mathbb{R})$. Since $ru\in
L^{1}_{\mathrm{loc}}(\mathbb{R})$, we see that $u^{\prime}\in L^{1}%
_{\mathrm{loc}}(\mathbb{R})$, hence $u\in W^{1,1}_{\mathrm{loc}}(\mathbb{R})$,
i.e. $u$ is absolutely continuous. But then $ru\in L^{2}_{\mathrm{loc}%
}(\mathbb{R})$, and $u\in H^{1}_{\mathrm{loc}}(\mathbb{R})$. Now we can
conclude that
\[
-(\partial_{x}+r)(\partial_{x}-r)u=(-\partial_{x}^{2} + q)u,
\]
where all operations should be applied in the distributional sense. We proved
that ${\mathfrak{D}}(H)\subset{\mathfrak{D}}(\widetilde{L}_{q})$ and $H$,
applied as a factorized operator (see \eqref{E:factor}), is a restriction of
$\widetilde{L}_{q}$.
\end{proof}

\begin{remark}
In particular, we can apply $H$ on $\mathfrak{D}(H)$ as $-\partial_{x}^{2}+q$
applied termwise, which is usually easier than to apply it in the factorized
form. To illustrate it, note that it may easily happen that $\mathfrak{D}%
(\widetilde{L}_{q})$ does not contain any function $u\in\mathcal{C}%
_{0}^{\infty}(\mathbb{R})$ (except $u\equiv0$). This is true e.g. for the
potential
\[
q(x)=\sum_{k=1}^{\infty}c_{k}\delta(x-x_{k}), \quad\sum_{k=1}^{\infty}%
c_{k}<\infty,
\]
where $c_{k}>0$ for all $k$, and the set $\{x_{k}\}_{k=1}^{\infty}$ is dense
in $\mathbb{R}$. Therefore, the same is true for $\mathfrak{D}(H)$.
\end{remark}

\begin{remark}
The operator $\widetilde{L}_{q}$ may be defined on the domain
\eqref{E:domain-exp} even without semi-boundednes requirement. But in this
case the resulting operator (called usually ``maximal operator") in
$L^{2}(\mathbb{R})$ will not necessarily be self-adjoint even if
$q\in\mathcal{C}^{\infty}(\mathbb{R})$ (see e.g. Sect. X.1 in \cite{RS:1975}
or Sect. II.1 and II.4.2 in \cite{BS:1991}).
\end{remark}

Now, note that $\lambda_{0}(q)$ (as defined in (\ref{eq.lambda})) is the
bottom of the spectrum of the self-adjoint operator $H$. So, if $\lambda
_{0}(q)>0$, then there exists a bounded, everywhere defined linear operator
$T:=H^{-1}$ in $L^{2}(\mathbb{R})$. It maps $L^{2}(\mathbb{R})$ onto
$\mathfrak{D}(H)$. We will now analyze the properties of the Schwartz kernel
of $T$, which is Green's function for the operator $H,$ where $H$ is an arbitrary, but fixed 
Schr\"odinger operator with a real potential $q \in H^{1}_{loc}({\mathbb R}^{n})$ such that $\lambda_0 (q) >0$.

\begin{lemma}
\label{L:estimates} The following estimates hold for any bounded open interval
$I\subset\mathbb{R}$ and any $u\in H^{1}_{\mathrm{loc}}(I)$:
\begin{equation}
\left\Vert u\right\Vert _{L^{\infty}(I)}\leq\left\vert I\right\vert
^{-1/2}\left\Vert u\right\Vert _{L^{2}(I)} +\left\Vert u^{\prime}\right\Vert
_{L^{1}(I)}, \label{eq.linfty}%
\end{equation}
and%
\begin{equation}
\left\Vert u^{\prime}\right\Vert _{L^{1}(I)}\leq\left\vert I\right\vert
^{1/2}\left\Vert Pu\right\Vert _{L^{2}(I)}+ \left\Vert r\right\Vert
_{L^{2}(I)}\left\Vert u\right\Vert _{L^{2}(I)}. \label{eq.Pinfty}%
\end{equation}
where $|I|$ means the length of the interval $I$.
\end{lemma}

\begin{proof}
Both estimates are well-known but we provide the proofs for the convenience
of the reader. We start with
\[
u(x)-u(y)=\int_{y}^{x}u^{\prime}(s)ds,\quad x,y\in I,
\]
which implies
\[
|u(x)-u(y)|\leq\Vert u^{\prime}\Vert_{L^{1}(I)},
\]
hence
\[
|u(x)|\leq|u(y)|+\Vert u^{\prime}\Vert_{L^{1}(I)},
\]
and
\[
\Vert u\Vert_{L^{\infty}(I)}\leq|u(y)|+\Vert u^{\prime}\Vert_{L^{1}(I)},
\]
for all $y\in I$. Integrating with respect to $y\in I$ and dividing by $I$, we
get
\[
\Vert u\Vert_{L^{\infty}(I)}\leq|I|^{-1}\Vert u\Vert_{L^{1}(I)}+\Vert
u^{\prime}\Vert_{L^{1}(I)}.
\]
Applying the Cauchy-Schwarz inequality in the first term in the right hand
side, we obtain \eqref{eq.linfty}.

From $u^{\prime}=Pu+ru$, taking $L^{1}$-norms of both sides and using the
Cauchy-Schwarz inequality we obtain
\[
\|u^{\prime}\|_{L^1(I)}\le\|Pu\|_{L^{1}(I)}+\|ru\|_{L^{1}(I)}\le
|I|^{1/2}\|Pu\|_{L^{2}(I)}+\|r\|_{L^{2}(I)}\|u\|_{L^{2}(I)},
\]
which proves \eqref{eq.Pinfty}.
\end{proof}

\begin{lemma}
\label{L:operest} Let us assume that $q \in H^{-1}_{loc}(\mathbb R)$ with $\lambda_{0}(q)>0$. Then

\noindent1) The operator $\overline{P}T$ is defined everywhere in
$L^{2}(\mathbb{R})$ and $\|\overline{P}T\| = \lambda_{0}(q)^{-1/2}$.

\noindent2) $\overline{P}TP^{*}=I$ on $\mathfrak{D}(P^{*})$.
\end{lemma}

\begin{proof}
1) Since $P^{\ast}(\overline{P}T)=(P^{\ast}\overline{P})T=I$, the operator
$\overline{P}T$ is everywhere defined in $L^{2}(\mathbb{R})$. Also, for any
$u\in L^{2}(\mathbb{R})$,
\[
\left\Vert \overline{P}Tu\right\Vert ^{2}=(\overline{P}Tu,\overline
{P}Tu)=(TP^{\ast}\overline{P}Tu,u)=(Tu,u)= \|T^{1/2}u\|^{2},
\]
so the first statement immediately follows. (In fact, the presentation
$\overline{P}T=UT^{1/2}$, with $U=\overline{P}T^{1/2}$, is the polar
decomposition of $\overline{P}T$.)

2) For any $u,v\in\mathfrak{D}(P^{*}\overline{P})$ we have
\[
((\overline{P}TP^{*})\overline{P}u,\overline{P}v) =(\overline{P}%
T(P^{*}\overline{P})u,\overline{P}v) =(\overline{P}u,\overline{P}v),
\]
so $\overline{P}TP^{*}u=u$ for all $u\in\mathfrak{D}(P^{*}\overline{P})$. It
remains to recall that $\mathfrak{D}(P^{*}\overline{P})$ is dense in
$L^{2}(\mathbb{R})$.
\end{proof}

Heuristically, Green's function $G=G(x,y)$ should be given by
\begin{equation}
G(x,y)=(T\delta(\cdot-y))(x).\label{E:Green}%
\end{equation}
So it is expected to satisfy
\begin{equation}
(-\partial_{x}^{2}+q(x))G(x,y)=\delta(x-y),\label{E:Green-eq}%
\end{equation}
and be the Schwartz kernel of a bounded linear operator in $L^{2}(\mathbb{R}%
)$. More precisely, we will prove:

\begin{lemma}
\label{lemma.green} Let us assume that $q \in H^{-1}_{loc}(\mathbb R)$ with $\lambda_{0}(q)>0$.
Then there exists a measurable, real-valued function $G(x,y)$ on
$\mathbb{R}\times\mathbb{R}$ with the following properties: \newline(i) For
any bounded interval $I$ in $\mathbb{R}$, there is a positive constant $C(I)$
so that
\[
\sup_{x\in I}\left(  \int G(x,y)^{2}~dy\right)  ^{1/2}\leq C(I)
\]
and the map $x\mapsto G(x,~\cdot~)$ is H\"{o}lder continuous of order $1/2$ as
a mapping from $I$ into $L^{2}(\mathbb{R})$. \newline(ii) For any $f\in
L^{2}(\mathbb{R})$,
\[
\left(  Tf\right)  (x)=\int G(x,y)~f(y)~dy
\]
(iii) $G(x,y)=G(y,x)$ almost everywhere. \newline(iv) For each fixed $x$, the
function $G(x,~\cdot~)$ belongs to $\mathfrak{D}(\overline{P})$. \newline(v)
For any $\varphi\in\mathcal{C}_{0}^{\infty}(\mathbb{R})$ and any
$x\in\mathbb{R},$ $\left(  \overline{P}G(x,\cdot),\overline{P}\varphi\right)
=\varphi(x).$
\end{lemma}

\begin{remark}
Part (v) states that $G(x,y)$, viewed as a function of $y$ with $x$ as
parameter, solves the equation $L_{q}G(x,y)=\delta_{x}(y)$ in distribution sense.
\end{remark}

\begin{proof}
[Proof of Lemma \ref{lemma.green}]In what follows, $r \in L^{2}_{loc}(\mathbb R)$ 
is fixed and satisfies $q = r' +r^{2}$, $I$ denotes a bounded
interval in $\mathbb{R},$ and $C(I)$ denotes a generic constant depending on
$\left\vert I\right\vert $ and $r$. Its value may vary from line to line. 

We will make repeated use of the following observation, based on Lemma
\ref{L:estimates}. If $\psi\in\mathcal{D}(\overline{P})$ and $I$ is a bounded
interval, then by Lemma \ref{L:estimates},
\begin{equation}
\sup_{x\in I}\left\vert \psi(x)\right\vert \leq C(I)\left(  \left\Vert
\psi\right\Vert _{L^{2}(I)}+\left\Vert \overline{P}\psi\right\Vert
_{L^{2}(\mathbb{R})}\right).  \label{eq.DP.sup}%
\end{equation}
 Using the boundedness of $\psi$ and the fact that
$\overline{P}\psi\in L^{2}(\mathbb{R})$, we can then deduce that
\begin{align}
\left\Vert \psi^{\prime}\right\Vert _{L^{2}(I)} &  \leq\left\Vert \overline
{P}\psi\right\Vert _{L^{2}(\mathbb{R})}+\left\Vert r\right\Vert _{L^{2}%
(I)}\left\Vert \psi\right\Vert _{L^{\infty}(I)}\label{eq.DP.Der}\\
&  \leq C(I)\left(  \left\Vert \psi\right\Vert _{L^{2}(\mathbb{R})}+\left\Vert
\overline{P}\psi\right\Vert _{L^{2}(\mathbb{R})}\right)  .\nonumber
\end{align}

Since, by Lemma \ref{L:operest}, $\left\Vert \overline{P}T\psi\right\Vert
_{L^{2}(\mathbb{R})}\leq C\left\Vert \psi\right\Vert _{L^{2}(\mathbb{R})}$, it
follows from (\ref{eq.DP.sup}) and (\ref{eq.DP.Der}) that for any $\psi\in
L^{2}(\mathbb{R})$, one has $T\psi\in H_{\mathrm{loc}}^{1}(\mathbb{R})$ with%
\begin{equation}
\sup_{x\in I}\left\vert (T\psi)(x)\right\vert \leq C(I)\left\Vert
\psi\right\Vert _{L^{2}(\mathbb{R})}\label{eq.Tpsi.sup}%
\end{equation}
and
\begin{equation}
\int_{I}\left\vert \frac{d}{dx}(T\psi)(x)\right\vert ^{2}~dx\leq
C(I)\left\Vert \psi\right\Vert _{L^{2}(\mathbb{R})}^{2}. \label{eq.Tpsi.Der}%
\end{equation}
In particular, for each $x$, the map
$x\mapsto\left(  T\psi\right)  (x)$ is a bounded linear functional on
$L^{2}(\mathbb{R})$. It follows from the Riesz representation theorem that
there is an element $G_{x}$ of $L^{2}(\mathbb{R})$ with
\[
(T\psi)(x)=(\psi,G_{x}).
\]
We claim that, also, the map $I\ni x\mapsto G_{x}\in L^{2}(\mathbb{R})$ is
H\"{o}lder continuous of order $1/2$. To see this, we use (\ref{eq.Tpsi.Der})
together with the Cauchy-Schwarz inequality to conclude that for $x$ and $y$
belonging to $I$,%
\[
\left\vert (T\psi)(x)-(T\psi)(y)\right\vert \leq C(I)\left\vert x-y\right\vert
^{1/2}\left\Vert \psi\right\Vert _{L^{2}(\mathbb{R})}%
\]
and thus%
\[
\sup\left\{  \left\vert \left(  \psi,G_{x}-G_{y}\right)  \right\vert :\psi\in
L^{2}(\mathbb{R}),~\left\Vert \psi\right\Vert _{L^{2}(\mathbb{R})}=1\right\}
\leq C(I)\left\vert x-y\right\vert ^{1/2}.
\]
This proves the required H\"{o}lder continuity. It follows that the map
$x\mapsto G_{x}$ is a weakly measurable map from $I$ into $L^{2}(\mathbb{R})$
with $\left\Vert G_{x}\right\Vert _{L^{2}(\mathbb{R})}$ bounded uniformly in
$x\in I$, so that $x\mapsto G_{x}$ may be regarded as an element of the space
$L^{2}(I;L^{2}(\mathbb{R}))$ consisting of weakly measurable,
square-integrable functions on $I$ taking values in $L^{2}(\mathbb{R})$. By
Theorem III.11.17 of \cite{DS:1958}, there is a measurable function
$G_{I}(x,y)$ on $I\times\mathbb{R}$ with the property that $G_{I}%
(x,~\cdot~)=G_{x}$ for every $x\in I$. As%
\[
(\varphi,T\psi)=\int_{I\times\mathbb{R}}\varphi(x)G_{I}(x,y)\psi(y)~dy~dx
\]
for any $\varphi\in L^{\infty}(I)$ and $\psi\in L^{2}(\mathbb{R})$, it is easy
to see that for any bounded intervals $I$ and $J$ with $I\subset J$, the
restriction of $G_{J}$ to $I\times\mathbb{R}$ equals $G_{I}$ almost everywhere
with respect to product measure on $I\times\mathbb{R}$. Taking a sequence of
bounded intervals $\left\{  I_{n}\right\}  $ with $I_{n}\nearrow\mathbb{R}$ as
$n\rightarrow\infty$, we can construct a measurable function $G$ on
$\mathbb{R}\times\mathbb{R}$ that obeys properties (i) and (ii). Property
(iii) follows from the symmetry of $T$.

To prove property (iv), let $\varphi\in\mathcal{D}(P^{\ast})$ and note that%
\[
(G_{x},P^{\ast}\varphi)=(TP^{\ast}\varphi)(x).
\]
By Lemma \ref{L:operest}, $\left\Vert \overline{P}TP^{\ast}\varphi\right\Vert
_{L^{2}(\mathbb{R})}\leq\left\Vert \varphi\right\Vert _{L^{2}(\mathbb{R})}$
holds. Hence $TP^{\ast}\varphi\in W_{\mathrm{loc}}^{1,1}%
(\mathbb{R})$ and, for any bounded interval $I$,%
\[
\sup_{x\in I}\left\vert \left(  TP^{\ast}\right)  (x)\right\vert \leq
C(I)\left\Vert \varphi\right\Vert _{L^{2}(\mathbb{R})}%
\]
so that%
\[
\left\vert \left(  G_{x},P^{\ast}\varphi\right)  \right\vert \leq
C(I)\left\Vert \varphi\right\Vert _{L^{2}(\mathbb{R})}%
\]
for any $\varphi\in\mathcal{D}(P^{\ast})$. This shows that $G_{x}%
\in\mathcal{D}(P^{\ast\ast})=\mathcal{D}(\overline{P})$, proving (iv).

Finally, to prove (v), let $\varphi\in\mathcal{D}(H)$ and compute%
\begin{align*}
\varphi(x) &  =(TH\varphi)(x)\\
&  =(H\varphi,G_{x})\\
&  =(\overline{P}\varphi,\overline{P}G_{x}).
\end{align*}
Since $\mathcal{D}(H)$ is dense in $\mathcal{D}(\overline{P})$ and point
evaluations are continuous in $\mathcal{D(}\overline{P})$, it follows that
$\varphi(x)=(\overline{P}\varphi,\overline{P}G_{x})$ for all $\varphi
\in\mathcal{D}(\overline{P})$. Since $\mathcal{C}_{0}^{\infty}(\mathbb{R}%
)\subset\mathcal{D}(\overline{P})$, (v) is proved.
\end{proof}

To construct positive solutions from Green's function, we will need the
following lemma.

\begin{lemma}
\label{lemma.zeros} Suppose $\lambda_{0}(q)>0$ and that $y\in H_{\mathrm{loc}%
}^{1}(\mathbb{R})$, $L_{q}y=0$, and either \newline(i) $y\in L^{2}(0,\infty)$
and $Py\in L^{2}(0,\infty)$, or \newline(ii) $y\in L^{2}(-\infty,0)$ and
$Py\in L^{2}(-\infty,0)$ . \newline Then, either $y$ has no zeros on
$\mathbb{R}$ or $y$ is identically zero on $\mathbb{R}$.
\end{lemma}

\begin{proof}
We will give the proof assuming (i) holds since the proof assuming (ii) holds
is similar. Suppose that $y\in H_{\mathrm{loc}}^{1}(\mathbb{R)}$ solves
$L_{q}y=0$, (i) holds, and $y(x_{0})=0$ for some $x_{0}\in\mathbb{R}$. We will
assume without loss that $x_{0}=0$. By assumption (i), the function
\[
w(x)=\left\{
\begin{array}
[c]{ll}%
y(x), & 0\leq x<\infty\\
0, & x<0
\end{array}
\right.
\]
belongs to $L^{2}(\mathbb{R})\cap H_{\mathrm{loc}}^{1}(\mathbb{R})$, and
$Pw\in L^{2}(\mathbb{R})$, hence $w\in W_{\mathrm{loc}}^{1,1}(\mathbb{R})$. It
follows from Lemma \ref{operatorP} that $w\in\mathfrak{D}(\overline{P}).$

We claim that there is a sequence $\left\{  \varphi_{n}\right\}  $ from
$\mathcal{C}_{0}^{\infty}(\mathbb R)$ with support contained in $(0, \infty)$
 so that $\varphi_{n}\rightarrow w$ and
$P\varphi_{n}\rightarrow\overline{P}w$ in $L^{2}(\mathbb{R})$. If so then, on
the one hand,%
\begin{equation}
\left(  P\varphi_{n},\overline{P}w\right)  =0 \label{eq.pre.Pw}%
\end{equation}
since $L_{q}y=0$ and $L_{q}y=L_{q}w$ as distributions on $(0,\infty).$
Taking limits in (\ref{eq.pre.Pw}) as $n\rightarrow\infty$
we have%
\begin{equation}
(\overline{P}w,\overline{P}w)=0, \label{eq.Pw}%
\end{equation}
hence $\overline{P}w=0$. On the other hand, we have%
\begin{equation}
(P\varphi_{n},P\varphi_{n})\geq\lambda_{0}(q)\left\Vert \varphi_{n}\right\Vert
^{2} \label{eq.Pphin}%
\end{equation}
where $\lambda_{0}(q)>0$. Taking limits as $n\rightarrow\infty$ in
(\ref{eq.Pphin}) and using (\ref{eq.Pw}), we conclude that $w=0$. It follows
from the uniqueness of solutions to $L_{q}y=0$ with prescribed initial data
that $y=0$ identically.

Thus, it remains to prove the existence of a sequence $\left\{  \varphi
_{n}\right\}  $ with the claimed properties. First, we show that $w$ may be
approximated by functions which vanish identically near $x=0$. Let $\chi
\in\mathcal{C}^{\infty}(\mathbb{R})$ with $0\leq\chi(x)\leq1$, $\chi(x)=0$ for
$x\leq1$, $\chi(x)=1$ for $x\geq2$, and $|\chi'(x) | \le 2$ for all $x \in \mathbb R$. Let $\chi_{\varepsilon}(x)=\chi
(x/\varepsilon)$. The functions $w_{\varepsilon}(x)=\chi_{\varepsilon}(x)w(x)$
converge to $w$ in $L^{2}(\mathbb{R})$ as $\varepsilon \to 0$ by dominated convergence. We claim
that, also, $Pw_{\varepsilon}\rightarrow \overline{P}w$ in $L^{2}(\mathbb{R})$. To see
this, compute%
\begin{equation}
Pw_{\varepsilon}=\chi_{\varepsilon}^{\prime}(x)w(x)+\chi_{\varepsilon
}(x)(\overline{P}w)(x). \label{eq.Pweps}%
\end{equation}
The second term in the right-hand side of (\ref{eq.Pweps}) converges to $\overline{P}w$ in
$L^{2}(\mathbb{R})$ by dominated convergence while the first one converges to $0$ in $L^{2}(\mathbb R)$ by the following reasons.
Observing that
\[
   \int\chi^{\prime}_{\varepsilon}(x)^{2}\left\vert w(x)\right\vert ^{2}dx
   \le   \frac{4}{\varepsilon^{2}}   \int_{\varepsilon}^{2\varepsilon}
   \left( \int_{0}^{x} \left\vert w'(t)\right\vert dt \right)^{2}dx
   \le 4 \int_{0}^{2\varepsilon} \left\vert w'(t)\right\vert^{2}dt,
\]
we conclude
\[
 \int\chi^{\prime}_{\varepsilon}(x)^{2}\left\vert w(x)\right\vert ^{2}dx
   \to 0 \, \text{ as } \varepsilon \to 0.
\]
\ Thus $Pw_{\varepsilon}\rightarrow\overline{P}w$ in $L^{2}(\mathbb{R})$.

Letting $\varepsilon=1/n$, the function $w_{1/n}$ has support in
$[1/n,\infty)$. We can use smooth cut-off functions and Friedrichs mollifiers
as in the proof of Lemma \ref{operatorP} to find a $\mathcal{C}_{0}^{\infty}$
function $\varphi_{n}$ with support in $[1/(2n),\infty)$ so that $\left\Vert
\varphi_{n}-w_{1/n}\right\Vert <1/n$ and $\left\Vert P\varphi_{n}-\overline
{P}w_{1/n}\right\Vert <1/n$. In this way we obtain a sequence $\left\{
\varphi_{n}\right\}  $ from $\mathcal{C}_{0}^{\infty}(0,\infty)$ so that
$\left\Vert \varphi_{n}-w\right\Vert \rightarrow0$ and $\left\Vert
P\varphi_{n}-\overline{P}w\right\Vert \rightarrow0$ as $n\rightarrow\infty$.
\end{proof}

As an application of the results obtained in this appendix, we give an alternative proof
of Proposition \ref{prop.2pos}

\textit{Proof of Proposition \ref{prop.2pos}.} We claim that there exists an
$x\in\mathbb{R}$ so that $y\mapsto$ $G(x,y)$ does not vanish identically on
$(x,\infty)$. If not then $G(x,y)=0$ for all $\left(  x,y\right)  $ with $y>x$
and hence, by Lemma \ref{lemma.green}(iii), for all $y\neq x$. Therefore
$G(x,y)=0$ a.e., a contradiction. Now choose such an $x$. Then the function
$\psi_{+}(y)=G(x,y)$ for $y>x$ is not identically zero on $(x,\infty)$. From
Lemma \ref{lemma.green}(i), (iv), and (v), $\psi_{+}(y)\in L^{2}\left(
x,\infty\right)  $, $P\psi_{+}\in L^{2}(x,\infty)$, and $L_{q}\psi_{+}=0$ for
$y>x$. Let $Q$ be an antiderivative of $q$ and let $\left\{  y_{+}%
,u_{+}\right\}  $ be the unique solution to the system (\ref{eq.sys}) with
initial data $y_{+}(x+1)=\psi_{+}(x+1)$ and $\left(  u_{+}\right)
(x+1)=\left(  \psi_{+}-Q\psi_{+}\right)  (x+1)$. Then $y_{+}$ coincides with
$\psi_{+}$ on $\left(  x,\infty\right)  $, so $y_{+}$ and $Py_{+}$ belong to
$L^{2}\left(  0,\infty\right)  $. It follows from Lemma \ref{lemma.zeros} that
$y_{+}$ has no zeros, so by changing signs if necessary we conclude that
$y_{+}\in L^{2}(0,\infty)$ and $y_{+}$ is strictly positive on $\mathbb{R}$. A
similar construction considering the function $\psi_{-}(y)=G(x,y)$ for some
$x\in\mathbb{R}$ and $y<x$ leads to a strictly positive solution $y_{-}$ of
$L_{q}y=0$ with $y_{-}\in L^{2}(-\infty,0)$. If $y_{+}$ and $y_{-}$ were
linearly dependent, then after multiplying one of them by an appropriate
constant, we would obtain a function $\psi$ in the domain of $H$, $\psi$ not
identically zero, with $H\psi=0$, which is impossible since $\lambda_{0}%
(q)>0$. Thus $y_{+}$ and $y_{-}$ are linearly independent.%
\hfill$\Box$%

\section{Related Work}

\label{app.work}

The Miura map was introduced by Miura \cite{Miura:1968}, \cite{Miura:1976} and
played an important role in the search of integrals of motion for the
Korteweg-de Vries equation. Miura discovered that his map takes smooth
solutions of mKdV to smooth solutions of KdV. Hence it can serve as a tool to
derive results on the initial value problem for KdV from results on the
initial value problem for mKdV -- see e.g. \cite{Tsu:1989}. Despite the fact
that the Miura map is not one-to-one, when considered, for example, as a map
between appropriate Sobolev spaces, it is also possible to use it to derive
results for the initial value problem of mKdV from results of the initial
value problem of KdV -- see e.g. \cite{GSS:1991}, \cite{CKSTT:2003},
\cite{KT:2003c}.
\newline\textit{{Miura map on the circle:}}\ Motivated by
earlier work of Ambrosetti and Prodi \cite{AP:1972} on certain nonlinear
elliptic boundary value problems, McKean and Scovel \cite{MS:1986} studied -
among other nonlinear maps - the Miura map on the circle $\mathbb{T}.$ They
exhibited a global fold structure for the Miura map when viewed as a map from
$H^{1}(\mathbb{T})$ to $L^{2}(\mathbb{T}).$ For further results in this
direction, see Bueno and Tomei \cite{BT:2002}. Later, Korotyaev
\cite{Korotyaev:2001}, \cite{Korotyaev:2003} and Kappeler and Topalov
\cite{KT:2003a} extended the global fold picture to the Miura map from
periodic functions in $L^{2}(\mathbb{T})$ to $H^{-1}(\mathbb{T})$. Kappeler
and Topalov proved existence and well-posedness of solutions to the mKdV
equation with initial data in $L^{2}(\mathbb{T})$ \cite{KT:2003c}, using
\cite{KT:2003a} and their results on the initial value problem for the
periodic KdV\ equation established in \cite{KT:2003bc}. \newline\textit{{Miura
map on the line:}}\ On the line, the Miura map and related topics have also
been investigated extensively, and not exclusively with a view towards
applications for solving the initial value problem of KdV or mKdV.
\newline$\bullet$ Positive solutions of Schr\"{o}dinger equations or more
generally of second order elliptic equations -- in particular in connection with
spectral properties of the corresponding operators -- have been extensively
studied in various settings. We only mention the result, referred to as
Allegretto-Piepenbrink theorem in \cite{CFKS:1987}, Theorem 2.12 or in
\cite{Simon:1982}, section C.8. This theorem states that for potentials $q\in
L_{\mathrm{loc}}^{1}(\mathbb{R}^{n}),$ satisfying some additional conditions,
$(-\Delta+q)u=\lambda u$ has a nonzero solution $u$ (in the sense that $u\in
W_{\mathrm{loc}}^{2,1}(\mathbb{R}^{n})$ and $qu\in L_{\mathrm{loc}}%
^{1}(\mathbb{R}^{n})$), which is nonnegative everywhere, if and only if
$\inf(\operatorname*{spec}$($-\Delta+q))\geq\lambda$. See \cite{CFKS:1987} or
\cite{Simon:1982} for further details and references to the papers of
Allegretto and of Piepenbrink as well as additional references. In the
one-dimensional case at hand, the equivalence of the statements (ii) and (iii)
of Theorem \ref{thm.im1} for potentials $q\in L_{\mathrm{loc}}^{1}%
(\mathbb{R})$ is well known -- see \cite{Hartman:1982}, Theorems XI.6.1 and
XI.6.2 and Corollary XI.6.1, \cite{Murata:1986}, Appendix 1, or \cite{GZ:1991}%
, Theorem 3.1. Thus Theorem \ref{thm.im1} as stated above shows in particular
that this equivalence continues to hold for $q\in H_{\mathrm{loc}}%
^{-1}(\mathbb{R})$. \newline$\bullet$ With regard to the characterization of
the image of the Miura map $B_{0}$ , we mention the result of Ablowitz et. al.
\cite{AKS:1979} which characterizes the image of Schwartz space by $B_{0}$ in
terms of the scattering data of these potentials as well as the result of
Tsutsumi \cite{Tsu:1989}, stating that any finite, positive Radon measure is
in the image of the Miura map $B_{0}:L^{2}(\mathbb{R})\rightarrow
H^{-1}(\mathbb{R}).$ Further, the case where $q$ is continuous is treated by
Hartman \cite{Hartman:1982}, Chapter XI.7, Lemma 7.1. The result stated in
Theorem \ref{thm.im2} sharpens all these results and puts them into a broader
perspective. \newline$\bullet$ The dichotomy described in Theorem
\ref{thm.im3} has another interpretation which does not involve the Miura map
at all: Murata \cite{Murata:1986}, Appendix 1, describes the dichotomy stated
in Theorem \ref{thm.im3} -- again for potentials in $q\in L_{\mathrm{loc}}%
^{1}(\mathbb{R})$ -- in terms of the notion of subcritical, critical, and
supercritical potentials, where in his terminology (i) $q$ is called
subcritical if $L_{q}$ has a positive Green's function, (ii) $q$ is called
critical if $L_{q}\geq0$ and does \textit{{not }}have a positive Green's
function , and (iii) $q$ is called supercritical if $L_{q}$ is \textit{{not }%
}nonnegative. See Simon \cite{Simon:1981} for an alternative notion of
subcritical and critical potentials. In \cite{Murata:1986}, Theorem A.5,
Murata shows that (i) $q\in L_{\mathrm{loc}}^{1}(\mathbb{R})$ is subcritical
iff $L_{q}y=0$ admits two linearly independent positive solutions in
$W_{\mathrm{loc}}^{2,1}(\mathbb{R})$ and that (ii) $q\in L_{\mathrm{loc}}%
^{1}(\mathbb{R})$ is critical iff $L_{q}y=0$ has up to scaling one positive
solution $y\in W_{\mathrm{loc}}^{2,1}(\mathbb{R}).$ These results of Murata
for one-dimensional Schr\"{o}dinger operators were later extended by Gesztesy
and Zhao \cite{GZ:1991}, Theorem 3.6, to more general Sturm-Liouville
operators. Our results on the dichotomy for Schr\"{o}dinger operators obtained
in this paper extend the results of Murata (and of Gesztesy and Zhao) in two
directions. First, we consider potentials which are real-valued distributions
in a Sobolev space with negative index of smoothness $\beta\geq-1$. Hence they
are not necessarily functions. Second, we describe
geometric aspects of the dichotomy: see Theorem \ref{thm.im3} and Theorem
\ref{thm.nofold}. 
\newline\textit{{Schr\"{o}dinger operators with singular
potentials:}}\ Recently, the operators $L_{q}$, considered on an interval
$(a,b)$, with potential $q$ in a Sobolev space with negative index of
smoothness, have been studied by various authors. In particular, we mention the
paper \cite{SS:2003} where different approaches to define the operator $L_{q}$
are discussed in detail and asymptotics for the eigenvalues and eigenfunctions
of these operators are obtained. See also \cite{HM:2001}, \cite{HM:2002},
\cite{KM:2001}, \cite{Korotyaev:2001}, \cite{Korotyaev:2003}, \cite{NS:1999},
\cite{SS:1999}, \cite{SS:2003} as well as \cite{SS:2003} for further
references. 
\newline\textit{{Initial value problem for KdV:}}\ The initial
value problem for KdV on the line has been extensively studied. We only
mention that, based on the works of Bourgain \cite{Bourgain:1993a},
\cite{Bourgain:1993b}, it has been proved by Kenig-Ponce-Vega \cite{KPV:1996}
that KdV is locally uniformly $C^{0}$ well-posed on $H^{s}(\mathbb{R})$ for
$s>-3/4$ and later, by Colliander, Keel, Staffilani, Takaoka and Tao
\cite{CKSTT:2003}, that KdV is \textit{{globally}} uniformly $C^{0}$
well-posed on $H^{s}(\mathbb{R})$ for $s>-3/4$. An existence result for the
limiting case $s=-3/4$ has been obtained by Christ, Colliander and Tao
\cite{CCT:2003}. Beside the work of Tsutsumi already mentioned above on
solutions of KdV with positive Radon measures as initial data, it has been
shown in \cite{Kappeler:1986} that for measures of bounded variation with
sufficient decay at infinity as initial data, \ there exists a
\emph{classical} solution for $t>0$.

\end{document}